\documentclass[a4paper,11pt]{article}
\usepackage[catalan,spanish,english]{babel}
\usepackage[spanish,english]{babel}
\usepackage[psamsfonts]{amssymb}
\usepackage{cmmib57}
\usepackage{exscale}
\usepackage{amsmath}
\usepackage{amscd}
\usepackage{latexsym}
\usepackage{graphicx}
\usepackage{amsmath}
\newtheorem{theorem}{Theorem}

\newtheorem{definition}[theorem]{Definition}
\newtheorem{example}[theorem]{Example}

\newtheorem{lemma}[theorem]{Lemma}

\newtheorem{proposition}[theorem]{Proposition}
\newtheorem{remark}[theorem]{Remark}

\newenvironment{proof}[1][Proof]{\textbf{#1.} }{\ \rule{0.5em}{0.5em}}

\newcommand{\cf}{\mathcal{F}}
\newcommand{\ch}{\mathcal{H}}

\newcommand{\cd}{\mathcal{D}}
\newcommand{\cc}{\mathcal{C}}
\newcommand{\cl}{\mathcal{L}}
\newcommand{\cb}{\mathcal{B}}

\newcommand{\cp}{\mathcal{P}}
\newcommand{\ct}{\mathcal{T}}
\newcommand{\ppa}{]\!]}
\newcommand{\ppc}{[\![}

\newcommand{\R}{\mathbb{R}}
\newcommand{\N}{\mathbb{N}}
\newcommand{\Z}{\mathbb{Z}}
\newcommand{\HH}{\mathbb{H}}

\newcommand{\lla}{\left\langle}
\newcommand{\rra}{\right\rangle}

\newcommand{\no}{\noindent}
\newcommand{\ds}{\displaystyle}

\begin{document}
\begin{titlepage}
\null
\vspace{2cm}
\begin{center}
{\LARGE \bf SPDEs with coloured noise:\\[2mm]
Analytic and stochastic approaches}\\[2mm]

\bigskip

by\\
\vspace{7mm}

\begin{tabular}{l@{\hspace{10mm}}l@{\hspace{10mm}}l}
{\sc Marco Ferrante}$\,^{(\ast)}$  &and&{\sc Marta Sanz-Sol\'e}$\,^{(\ast\ast)}$\\
{\small Dipartimento di Matematica}&&{\small Facultat de Matem\`atiques}\\
{\small Universit\`a di Padova}  &&{\small Universitat de Barcelona}\\
{\small Via Belzoni 7}  &&{\small Gran Via 585}\\ 
{\small 35131 Padova, Italy}  &&{\small 08007 BARCELONA, Spain}\\  
{\small e-mail:ferrante@math.unipd.it}  &&{\small e-mail: marta.sanz@ub.edu}\\            
\end{tabular}
\end{center}

\vspace{1.5cm}

{\bf Abstract:} We study strictly 
parabolic stochastic partial differential equations on $\R^d$, $d\ge 1$, driven
by a Gaussian noise white in time and coloured in space. Assuming that the coefficients
of the differential operator are random, we give sufficient conditions on the correlation
of the noise ensuring H\"older continuity for the trajectories of the solution of the equation.
For self-adjoint operators with deterministic coefficients, the mild and weak formulation
of the equation are related, deriving path properties of the solution to a parabolic Cauchy
problem in evolution form.
\vfill

\medskip

\medskip


\vspace{1.5cm}

\footnotesize
{\begin{itemize}
\item[$^{(\ast)}$] Partially supported by the grant COFIN 2001-015341-006 of MIUR, Italy. 
\item[$^{(\ast\ast)}$] Partially supported by the grant BFM2003-01345 from the \textit{Direcci\'on
General de Investigaci\'on, Ministerio de Ciencia y Tecnolog\'{\i}a}, Spain.
\end{itemize}}
\end{titlepage}

\section{ Introduction}

Stochastic partial differential equations (SPDEs) can be analized by different
approaches related with the classical deterministic  methods.
Let us mention the variational point of view (\cite{krylovrozovsky}, \cite{pardoux},
\cite{rozovsky}) and the semigroup approach (\cite{dpz}), based on analytical methods,
and the more genuine probabilistic
setting  using stochastic integration with respect to martingale measures  (\cite{w}, \cite{d1}).

The variational approach leads in particular to a now very complete $L_2$-theory (see \cite{rozovsky}).
However, this theory does not provide with sharp results on the properties of the trajectories of the solutions
of SPDEs, except  in the time variable. 
A more deep analytical insight into parabolic SPDEs  has been recently given by Krylov and
Lototsky, developing an $L_p$-theory with $p\in[2,\infty)$ (see \cite{k1}, \cite{k2} and the references herein,
\cite{kl}). 
This theory allows to obtain properties of the trajectories -both in time and space- quite sharp, using
Sobolev type imbeddings. Let us point out that in \cite{k1}, \cite{k2} the coefficients of the differential
operator can be random, therefore the theory applies to a very general class of equations. In a similar
spirit, parabolic SPDEs with deterministic coefficients in H\"older classes have been studied in \cite{mik}.

In this paper we study stochastic partial differential equations in the whole space $\R^d$,
with arbitrary dimension $d\ge1$,
driven by a Gaussian noise white in time and with homogeneous spatial correlation.
The differential operator is strictly parabolic with random coefficients, the
free terms are random as well. 
Using the analytical approach of \cite{k2} (see also \cite{k1}), we give sufficient
conditions on the correlation of the noise ensuring the existence of a solution with values
on some subspace of $L_p(\R^d)$, $p\in[2,\infty)$, and then, by means of Sobolev type imbeddings, we obtain
the existence of a random field, indexed by time and space, which is a version of the solution
and has its trajectories jointly H\"older continuous in $t, x$. 

A similar question using the evolution approach has been addressed in some previous articles.
In fact, parabolic equations with random coefficients in spatial dimension $d=1$, driven by  a space-time white noise 
have been studied in \cite{aln1}. The main result is the existence
of a {\it continuous} random field solution to the equation. 
The mild form of the equation contains a stochastic convolution with an anticipating
integrand. Therefore, the analysis requires tools of anticipating stochastic calculus
-a intricate machinery based on Malliavin calculus- and needs
 a strong regularity in terms of the random component $\omega$. These type
of hypothesis can be avoided with the analytic approach. In fact,  in this situation 
stronger results are given in \cite{k2}, Theorem 8.5 and Remark 8.7, where joint {\it H\"older continuity} is
obtained.

For $d\ge 1$ and driving noise of the same kind that the one we are considering
in this paper, joint H\"older continuity
for the stochastic heat equation  in its mild form has been  obtained in \cite{ss2}. The result, proved by means
of Kolmogorov's continuity criterium, is an
extension of the one stated in \cite{w} for $d=1$. 

The analytic approach considers the formal SPDE in a {\it weak} form (see (\ref{1.3}), (\ref{1.3.80})).
Studying the relationship between the {\it weak} and the {\it mild} formulation of the SPDE (see (\ref{2.4}))
 gives the possibility of transferring results obtained in the analytic setting to the evolution
scenary. The last part of the article is devoted to this topic, in the particular case where
the differential operator is self-adjoint and its
coefficients  are deterministic. In the framework of a $L_2$-theory, for a 
Neumann boundary-value problem with a strictly
parabolic divergence operator, this question has been studied in \cite{sv1}.

\section{Some preliminaries and notation}

We denote by $\cd(\R^{d+1})$ the space of Schwartz test functions
(\cite{schwartz}, page 24).
On a complete probability space $(\Omega,\cf,P)$, we consider
a Gaussian process $\{F(\phi),\phi\in \cd(\R^{d+1})\}$, mean zero,
with covariance functional given by
\begin{equation}
\label{1.1}
\begin{array}{l}
{\ds E(F(\phi),F(\psi))= \int_{\R_+}ds \int_{\R^d} \Gamma(dx)
(\phi(s,\cdot)\ast\widetilde{\psi}(s,\cdot))(x)}
\\ \\
{\ds = \int_{\R_+}ds \int_{\R^d} \mu(d\xi) \cf\phi(s,\cdot)(\xi)
\overline{\cf\psi(s,\cdot)(\xi)} .}
\end{array}
\end{equation}
In (\ref{1.1}), 
$\Gamma$ is a
non-negative, non-negative definite, tempered measure,
$\widetilde{\psi}(s,x) =\psi(s,-x)$, $\mu$ is the non-negative tempered measure on $\R^d$
defined by $\cf^{-1} \Gamma$, where $\cf$ denotes the 
Fourier transform operator. We notice that $\Gamma$ is a symmetric measure (\cite{schwartz}, Chap. VII,
Theorem XVII). 

For any test function $f,g\in \cd(\R^{d+1})$, the functional
\begin{equation*}
Q(f,g)=
\int_{\R^d}
\Gamma(dx) (f\ast\widetilde{g})(x)
\end{equation*}
is non-negative and
translation invariant, that means,
$Q(f,g)=Q(\tau_xf,\tau_xg)$,
where
$\tau_xf(\cdot)=f(\cdot+x)$
(see Gel'fand and Vilenkin \cite{gv1}, pag. 169).

Following Dalang and Frangos \cite{df1} (see also Dalang
\cite{d1}) the process $F$ can be extended to a worthy martingale
measure in the sense of Walsh. We will denote by $\{F(t,A),
t\geq0, A\in \cb_b(\R^d)\}$ this extension and by $\cf_t$ the
$\sigma$--field generated by $\{F(s,A),0\leq s\leq t, A\in
\cb_b(\R^d)\}$.

Consider the inner product on
$\cd(\R^{d})$ defined by
\[
\lla f,g \rra_\ch = \int_{\R^d} \Gamma(dx)
(f\ast\widetilde{g})(x).
\]
Let $\ch$ be the completion of $\cd(\R^{d+1})$ with respect to the
norm derived from $\lla\cdot,\cdot\rra_\ch$. For any complete
orthonormal system (CONS) $\{e_j,j\geq0\}\subset\cd(\R^{d+1})$ of $\ch$,  
define
\begin{equation}
\label{1.1.1}
W^k(t)=
\int_0^t
\int_{\R^d}
F(ds,dx)
e_k(x) ,
\end{equation}
$k\geq0$, where the integral must be understood in Walsh's sense. The process
$\{W^k(t), t\in[0,T], k\ge 0\}$ is  a sequence of
independent standard Brownian motions.

One can check that for any predictable
process $X$,
\begin{equation}
\label{1.12} \int_0^t \int_{\R^d} F(ds,dx) X(s,x) =
\sum_{k=0}^\infty \int_0^t W^k(ds) \lla X(s,\cdot),e_k(\cdot)
\rra_\ch .
\end{equation}
In particular, for any
$\phi\in \cd(\R^{d+1})$
\begin{equation}
\label{1.13} F(t,\phi):= \int_0^t \int_{\R^d} F(ds,dx) \phi(x) =
\sum_{k=0}^\infty \lla \phi,e_k \rra_\ch W^k(t) .
\end{equation}

Let $p\in (1,+\infty)$, $n\in \R$ and $d\in \N$. We denote by
$H^n_p=H^n_p(\R^d)$ the fractional Sobolev space
consisting of distributions $g$ on $\R^d$ such that there exists
$f\in L_p(\R^d)$ and $g=(1-\Delta)^{-\frac{n}{2}}f$. It is 
a Banach space
endowed with the norm
\[
\|u\|_{n,p}=\|(1-\Delta)^{n/2}u\|_p,
\]
where $\|\cdot\|_p$ denotes the usual norm  of $L_p(\R^d)$  and $\Delta$ 
is the Laplacian operator on $\R^d$. It is important to notice
that $\|\cdot\|_{n,p}\leq \|\cdot\|_{m,p}$ for $n\leq m$; this
gives rise to the embeddings
\[
\cdots \subset H^m_p \subset H^n_p \subset \cdots
\subset L_p \subset \cdots \subset
H^{-n}_p \subset H^{-m}_p \subset \cdots
\]
When $n\in \Z_+$, the spaces $H^n_p$ coincide with the classical
Sobolev spaces $W^n_p$. Moreover, the space $\cc_0^\infty$ of
infinitely differentiable functions with compact support is dense
in each $H^n_p$. 
We refer the reader to \cite{ada} and \cite{t1}
for an extensive account on these spaces.

\section{SPDEs with random coefficients}

In this section, we analyze a parabolic spde, with Lipschitz coefficients,
driven by a noise $F$
as has been described in Section 2, 
under the prespective of the general theory developed in \cite{k1},
\cite{k2}. More precisely, we exhibit a relationship between the covariance
measure $\Gamma$ and a fractional differentiability degree $\eta$ leading, a.s.,
to jointly continuous solutions in time and in space. 

The results might be considered as a complement of those in Section 8.3 in \cite{k2},
where the spatial dimension is $d=1$ and the driving noise, white in time and in space.
Their proof consists in showing that the assumptions of 
Theorem 5.1 in \cite{k2} (see also Theorem 3.6  in \cite{k1})
are satisfied.

For the sake of completeness, we start by quoting some basic material from \cite{k1},
\cite{k2}.

Consider a fractional Sobolev space $H_p^n$, with fixed $p\in[2,\infty)$, $n\in\R$.
For any $u\in H_p^n$, $\phi\in \cc_0^\infty$, we define
\begin{equation}
\label{1.2} (u,\phi)= \int_{\R^d} [(1-\Delta)^{n/2}u](x)
[(1-\Delta)^{-n/2}\phi](x) dx.
\end{equation}
Let $\tau$ be a
stopping time with respect to $(\cf_t)_{t\geq 0}$ and $\cp$ be the
predictable $\sigma$--field. Set
$\HH^n_p(\tau)=L_p(\ppa0,\tau \ppa,\cp,H^n_p)$,
$\HH^n_p:=\HH^n_p(\infty)$. The spaces $\HH^n_p(\tau)$ are
a kind of stochastic fractional Sobolev spaces.

 We also introduce the following
notation:
\begin{quote}
\textit{$(f,g)\in \cf_p^n(\tau)$ if and only if $f\in
\HH^n_p(\tau)$, $g\in \HH^{n+1}_p(\tau,l^2)$, and we set
$\|(f,g)\|_{\cf_p^n(\tau)}=
\|f\|_{\HH^n_p(\tau)}+\|g\|_{\HH^{n+1}_p(\tau,l^2)}$}
\end{quote}
where $\HH^{n+1}_p(\tau,l^2)$ correspond to the space
of square summable sequences of elements of $\HH^{n+1}_p(\tau)$.
We denote by $(w_k(t), t\in[0,T], k\ge 0)$ a sequence of independent standard
Wiener processes.

\begin{definition}
( Definition 3.1, \cite{k2})
\label{def1}
For a distribution valued function
$u\in \cap_{T>0} \HH^n_p(\tau \wedge T)$,
we write
$u\in \ch_p^n(\tau)$
if
\begin{equation*}
u_{xx}\in \HH^{n-2}_p(\tau),
u(0,\cdot)\in L_p(\Omega, \cf_0, H_p^{n-2/p})
\end{equation*}
and there exists
$(f,g)\in \cf_p^{n-2}(\tau)$ such that, for any
$\phi\in \cc_0^\infty$,
the equality
\[
(u(t,\cdot),\phi)=
(u(0,\cdot),\phi)+\int_0^t ds (f(s,\cdot),\phi) +
\sum_{k=1}^\infty \int_0^t
w^k(ds)
(g^k(s,\cdot),\phi)
\]
holds for all $\tau$ a.s.

We set
\[
\|u\|_{\ch_p^n(\tau)}=\|u_{xx}\|_{\HH^{n-2}_p(\tau)}+
\|(f,g)\|_{\cf_p^{n-2}(\tau)}+
(E\|u(0,\cdot)\|_{n-2/p,p}^p)^{1/p} .
\]
\end{definition}

Let us recall the result
 on existence and
uniqueness of solution for stochastic partial differential
equations of parabolic type driven by a sequence of independent
Wiener processes. First, we  introduce some notation, then
the assumptions and finally, the statement.

Fix $n\in \R$ and $\gamma\in [0,1[$ be such that $\gamma=0$ if
$n=0, \pm1, \pm2, \ldots$; otherwise, $\gamma>0$ and is such that
$|n|+\gamma$ is not an integer. Define
\[
B^{|n|+\gamma}=
\left\{
\begin{array}{ll}
B(\R^d) & \mbox{ if } n=0
\\
\cc^{|n|-1,1}(\R^d) & \mbox{ if } n=\pm1, \pm2, \ldots
\\
\cc^{|n|+\gamma}(\R^d) & \mbox{ otherwise} \quad,
\end{array}
\right.
\]
where $B(\R^d)$ is the Banach space of bounded functions on
$\R^d$, $\cc^{|n|-1,1}(\R^d)$ is the Banach space of $|n|-1$ times
continuously differentiable functions whose derivatives of
$(|n|-1)$--st order are Lipschitz; $\cc^{|n|+\gamma}(\R^d)$ are
H\"older spaces. The spaces $B^{|n|+\gamma}(l_2)$ are defined in
the obvious way.

Consider the following equation on
$\ppa 0,\tau\ppa$:
\begin{equation}
\label{1.3}
\begin{array}{l}
du(t,x)=\left[
a^{i,j}(t,x) u_{x^i,x^j} (t,x) +
f(t,x,u)\right] dt
\\ \\
+g^k(t,x,u) dw_t^k .
\end{array}
\end{equation}
Notice that, in comparison with equation (5.1) in Krylov \cite{k2}, we 
take here $\sigma^{ik} \equiv 0$.

By a solution to the Cauchy problem for equation (\ref{1.3})
with initial condition $u_0$, we mean a stochastic
process $u\in \ch_p^{n+2}(\tau)$ such that for any test function
$\phi\in \cc_0^\infty$,
\begin{align}
(u(t,\cdot),\phi) & = (u(0,\cdot),\phi)+ \int_0^t ds (a^{i,j}(s,\cdot) u_{x^i,x^j}(s,\cdot) + f(s,\cdot,u),\phi)\nonumber\\
& + \int_0^t
w^k(ds)
(g^k(s,\cdot,u),\phi), \label{1.3.80},
\end{align}
for all $t\in \ppc 0,\tau\ppa$.

Assume the following conditions on the differential operator and on the coefficients of
the equation:
\begin{description}
\item[(A1)]
For any
$i,j=1,\ldots,n$,
\[
a^{i,j} : \Omega\times\R_+\times\R^d \longrightarrow\R
\]
is $\cp\otimes \cb(\R^d)$-- measurable.
For any $\omega\in \Omega$ a.s. and $t\geq 0$,
we have $a^{i,j}(t,\cdot)\in B^{|n|+\gamma}$
and $\|a^{i,j}(t,\cdot)\|_{B^{|n|+\gamma}}\leq K$.
Moreover, there exist $K, \delta>0$,
such that for any $\omega\in \Omega$, $t\geq 0$,
$x,\lambda\in \R^d$, 
\[
\delta|\lambda|^2 \leq a^{i,j}(t,x) \lambda^i \lambda^j \leq
K |\lambda|^2 .
\]
\item[(A2)]
For any $u\in H_p^{n+2}$,
$f(t,\cdot, u)$, $g(t,\cdot, u)$
are predictable processes taking values in $H_p^n$
and $H_p^{n+1}(l_2)$, respectively.

\no
In addition,
\begin{enumerate}
\item
$(f(\cdot,\ast,0),g(\cdot,\ast,0))\in \cf_p^n(\tau)$,
\item
$f,g$ are a.s. continuous in the third variable $u$;
\item
for any $\varepsilon>0$, there exists $K_\varepsilon$ such that
for any $u,v\in H_p^{n+2}$, $t\ge 0$, 
\[
\|f(t,\cdot,u)-f(t,\cdot,v\|_{n,p} +
\|g(t,\cdot,u)-g(t,\cdot,v)\|_{n+1,p}
\]
\[
\leq
\varepsilon \|u-v\|_{n+2,p}+K_\varepsilon
\|u-v\|_{n,p},
\]
a.s.
\end{enumerate}
\end{description}

The next result is a particular version of Theorem 5.1 in 
\cite{k2}
\begin{theorem}
\label{th1}
Let $u_0\in L_p(\Omega,\cf_0,H_p^{n+2-2/p})$. Then
the Cauchy problem (\ref{1.3}) on $\ppa 0,\tau\ppa$ with initial
condition $u(0,\cdot)=u_0$ has a unique solution $u\in
\ch_p^{n+2}(\tau)$. This solution satisfies
\[
\|u\|_{\ch_p^{n+2}(\tau)} \leq
N
\left\{
\|f(\cdot,\ast,0)\|_{\HH_p^n(\tau)} +
\|g(\cdot,\ast,0)\|_{\HH_p^{n+1}(\tau,l_2)}
\right.
\]
\[
\left.
+
(E \|u_0\|^p_{n+2-2/p,p})^{1/p}
\right\} ,
\]
where the constant $N$ depends only on $d, n, \gamma, p, \delta, K, T$
and the function $K_\varepsilon$.
\end{theorem}


Consider now the equation
\begin{align}
du(t,x)& =\big[a^{i,j}(t,x) u_{x^i,x^j} (t,x)+ b^i(t,x) u_{x^i} (t,x)\nonumber\\
&+ f(t,x,u(t,x))\big] dt + h(t,x,u(t,x)) F(dt,x),\label{1.4}
\end{align}
with initial condition 
$u(0,x)=u_0(x)$,
where $t\in \R_+$, $x\in \R^d$ and $F$ is the Gaussian process
introduced in the preceding section. The coefficients $f, h$ are
random real functions defined on $\ppa0,\tau \ppa\times \R^d\times \R$.
Under suitable assumptions, we
shall prove that this equation can be set in the framework of
Theorem \ref {th1} and deduce H\"older continuity of the trajectories of
its unique solution.

Let us write (\ref{1.4}) into the form (\ref{1.3}). We consider a 
CONS $\{e_j, j\geq 0\}$ of $\ch$. We have,
\begin{align}
\lla h(t,\cdot,u),e_k \rra_\ch& = \int_{\R^d} \Gamma(dx)
(h(t,\cdot,u)\ast\tilde{e}_k)(x)\nonumber \\
&= \int_{\R^d} dy \ h(t,y,u)
\int_{\R^d} \Gamma(dx) \tilde{e}_k(x-y), \label{1.4.1}
\end{align}
where in the last equality we have applied
Fubini's theorem.

Therefore, the term $h(t,x,u(t,x)) F(dt,x)$ can be rewritten as
$g^k(t,x,u(t,x))$ $W^k(dt)$ with
\[
g^k(t,x,u(t,x))=
h(t,x,u(t,x))
\int_{\R^d} \Gamma(dy) \tilde{e}_k(y-x),
\]
and $W^k$ defined in (\ref{1.1.1}) (see (\ref{1.12})).

 Indeed, in the integral formulation,
the contribution of the last term in (\ref{1.4}) is,
for $\phi\in \cc_0^\infty$, $\int_0^t \int_{\R^d} F(dt,dx)\phi(x) h(t,x,u(t,x))$.
By virtue of  (\ref{1.12}) and (\ref{1.4.1}), 
\begin{align*}
&\int_0^t \int_{\R^d} F(dt,dx) \phi(x) h(t,x,u(t,x))  = 
 \sum_{k=0}^\infty \int_0^t W^k(ds) \lla h(t,\cdot,u(t,\cdot))\phi,
e_k\rra_\ch \\
& =\sum_{k=0}^\infty \int_0^t W^k(ds) \int_{\R^d}
h(t,y,u(t,y))\phi(y) \left(\int_{\R^d} \Gamma(dx)
\tilde{e}_k(x-y)\right)\\
& =\sum_{k=0}^\infty \int_0^t W^k(ds)
\left(h(t,\cdot,u(t,\cdot))\int_{\R^d} \Gamma(dx)
\tilde{e}_k(x-\cdot),\phi\right)_2,
\end{align*}
where $(\cdot,\cdot)_2$ denotes the inner product in $L^2(\R^d)$.

Set
\[
v_k(x)=
\int_{\R^d}
\Gamma(dy)
\tilde{e}_k(y-x) .
\]
The following lemma provides a useful tool 
to apply Theorem \ref{th1} to equation (\ref{1.4}), where
$h(t,x,u(t,x)) F(dt,x)$ is replaced by $g^k(t,x,u(x,t)) W^k(dt)$.

For any $\eta>0$, we denote by $R_{\eta,d}(x)$ the kernel of the
operator $(1-\Delta)^{-\eta/2}$ on $\R^d$, that is,
\[
\left[
(1-\Delta)^{-\eta/2}u\right](x)=
 R_{\eta,d}*u .
\]
It is well known that
\[
R_{\eta,d}(x)=C_{\eta,d} |x|^{\frac{\eta-d}{2}}
K_{\frac{d-\eta}{2}}(|x|) ,
\]
where $C_{\eta,d}$ is the reciprocal of $\pi^{d/2} 2^{(d+\eta-2)/2}\Gamma(\frac{\eta}{2})$
and $K_\nu$ is the modified Bessel function of the third kind
(see \cite{do1}).
Notice that $R_{\eta,d}(x)$ is a radial function and that
\[
R_{\eta_1,d} \ast R_{\eta_2,d} = R_{\eta_1+\eta_2,d}
\]
for any $\eta_1, \eta_2>0$. Hence,
\begin{equation}
\label{1.44} \nu_{\eta,d}:= \|R_{\eta,d}\|^2_{\ch} = \int_{\R^d}
\Gamma(dx) R_{2 \eta,d}(x).
\end{equation}

\begin{lemma}
\label{lm1}
Let $\eta\in (0,\infty)$, $d\in \N$ be such that
\begin{equation}
\label{1.5} \nu_{\eta,d} = \|R_{\eta,d}\|^2_{\ch} < \infty .
\end{equation}
Let $h\in L_p(\R^d)$,
$g^k=v_kh$.
Then $g=\{g_k,k\geq 0\} \in H_p^{-\eta}(l_2)$ and
\begin{equation}
\label{1.6}
\|g\|_{-\eta,p}=\|\overline{h}\|_p\leq C \|h\|_p ,
\end{equation}
with
\[
\overline{h}(x)=\|R_{\eta,d}(x-\cdot)h\|_\ch
\]
and $C=\nu_{\eta,d}^{1/2}$.
\end{lemma}

\begin{proof}
Fubini's theorem and Parseval's identity yield
\[
\begin{array}{l}
{\ds \|(1-\Delta)^{-\eta/2} g(x)\|^2_{l_2} = \sum_{k=0}^\infty
\left( (1-\Delta)^{-\eta/2} g^k(x)\right)^2}
\\
{\ds = \sum_{k=0}^\infty \left( (R_{\eta,d} \ast (v_k
h))(x)\right)^2}
\\
{\ds = \sum_{k=0}^\infty \left( \int_{\R^d}  dy  R_{\eta,d}(x-y) \left(
\int_{\R^d} \Gamma(dz) \tilde{e}_k(z-y) \right) h(y) \right)^2}
\\
{\ds = \sum_{k=0}^\infty \left( \int_{\R^d}  \Gamma(dz) \left( \int_{\R^d} dy
R_{\eta,d}(x-y) h(y) \tilde{e}_k(z-y)  \right)
\right)^2}
\\
{\ds = \sum_{k=0}^\infty \left( \int_{\R^d} \Gamma(dz) 
\left(R_{\eta,d}(x-\cdot) h \ast \tilde{e}_k \right)(z) 
\right)^2}
\\
{\ds = \sum_{k=0}^\infty \lla R_{\eta,d}(x-\cdot) h,{e}_k
\rra^2_{\ch} = \| R_{\eta,d}(x-\cdot) h \|^2_{\ch} .}
\end{array}
\]
Therefore,
\[
\begin{array}{ll}
\|g\|_{-\eta,p}
&
{\ds = \|(1-\Delta)^{-\eta/2} g\|_{L_p(l_2)} =
\left( \int_{\R^d} dx \|(1-\Delta)^{-\eta/2} g(x)\|^p_{l_2}
\right)^{1/p}}
\\
&
{\ds = \left( \int_{\R^d}  dx \|R_{\eta,d}(x-\cdot) h\|^p_\ch
\right)^{1/p} = \|\overline{h}\|_p .}
\end{array}
\]

The second part of (\ref{1.6}) is a consequence of H\"older's inequality.
Indeed, first we notice that, since $\Gamma$ is translation invariant,
\begin{equation}
\label{1.6.1}
\nu_{\eta,d} = ||R_{\eta,d}||^2_ {\ch} = ||R_{\eta,d}(x-\cdot)||^2_ {\ch}.
\end{equation}

Then, 
\begin{align*}
&\|\overline{h}\|_{p}^p
 = \int_{\R^d} dx \|R_{\eta,d}(x-\cdot) h\|^p_\ch\\
 &= \int_{\R^d} dx \left( \int_{\R^d}  \Gamma(dy)
\int_{\R^d} dz R_{\eta,d}(x-(y-z)) h(y-z)
\widetilde{R}_{\eta,d}(x-z) \widetilde{h}(z) \right)^{\frac{p}{2}}\\
&\leq \int_{\R^d} dx ( ||R_{\eta,d}(x-\cdot)||_ {\ch}^{\frac{p}{2}-1})\\
&\quad \times \left[ \int_{\R^d}  \Gamma(dy) \int_{\R^d}
dz R_{\eta,d}(x-(y-z)) \widetilde{R}_{\eta,d}(x-z) |h(y-z)
\widetilde{h}(z)|^{\frac{p}{2}} \right] .
\end{align*}

Thus, (\ref{1.6.1}), Fubini's theorem and Schwarz's inequality and the invariance of
Lebesgue measure imply
\[
\begin{array}{rl}
\|\overline{h}\|_{p}^p \leq & {\ds \nu_{\eta,d}^{\frac{p}{2}-1}
\int_{\R^d} dx \int_{\R^d} dz \int_{\R^d} \Gamma(dy)
R_{\eta,d}(y-z) R_{\eta,d}(z)}
\\
& \quad {\ds \times |h(y-z+x)|^{\frac{p}{2}}
|\widetilde{h}(z+x)|^{\frac{p}{2}}}
\\
\leq & {\ds \nu_{\eta,d}^{\frac{p}{2}-1} \int_{\R^d} \Gamma(dy)
\left( R_{\eta,d} \ast \widetilde{R}_{\eta,d} \right)(y) \Big(
\int_{\R^d} dx |h(y-z+x)|^{p} \Big)^{\frac{1}{2}}}
\\
& \quad {\ds \times \left( \int_{\R^d} dx |\widetilde{h}(z+x)|^{p}
\right)^{\frac{1}{2}}}
\\
=
&
\nu_{\eta,d}^{\frac{p}{2}}
\|h\|_{p}^p .
\end{array}
\]
This completes the proof of the lemma.
\end{proof}

\begin{remark}
\label{rem3} Let $ \Gamma(dx)=\delta_{\{0\}}(x)$ and thus,
$\|\cdot\|_{\ch} = \|\cdot\|_{2}$. In this particular case  (\ref{1.6}) has been
obtained in Lemma 8.4 of Krylov \cite{k2}.
\end{remark}

Proposition 4.4.1 in \cite{lev} establishes that, if 
\[
\int_{\R^d}
\frac{\mu(d \xi)}{(1+|\xi|^2)^\eta}
< + \infty .
\]
then (\ref{1.5}) holds true.

The behavior of the Bessel function $K_\nu$ is well-known (see for
instance Donoghue \cite{do1}). In fact, in  a neighborhood  $O^+$ of
$0$,
\[
K_\nu(r) \sim
\left\{
\begin{array}{ll}
\log(r) \ , \
&
\mbox{if} \ \nu=0,
\\ \\
r^{-|\nu|} \ , \
&
\mbox{if} \ \nu \neq 0.
\end{array}
\right.
\]
While away from zero,
\[
K_\nu(r)=C_\nu \ e^{-r} .
\]
This leads to the following conclusions, which have already appeared 
in previous discussions on different classes of spde's (for instance, in
\cite{ss1}, \cite{lev}).
\begin{enumerate}
\item
Assume $ 0<\eta<\frac{d}{2}$. Then
\[
\nu_{\eta,d}< +\infty \Leftrightarrow
\int_{O^+} |x|^{2\eta-d} \Gamma(dx) < +\infty .
\]
\item
Let $ \eta=\frac{d}{2}$. Then
\[
\nu_{\eta,d}< +\infty \Leftrightarrow
\int_{O^+} |x|^{\frac{2\eta-d}{2}}
\log \left( \frac{1}{|x|} \right)
\Gamma(dx) < +\infty .
\]
\item
If $ \eta>\frac{d}{2}$. Then
$ \nu_{\eta,d}< +\infty$, without any additional condition on $\Gamma$.
\end{enumerate}

\begin{example}
{\bf (Riesz kernels)}: Set $ \Gamma(dx)=|x|^{-\alpha} dx$, with
$\alpha\in(0,d)$. Then, for $\eta \in (0,\frac{d}{2}]$,
$\nu_{\eta,d} < +\infty$ if and only if $ \alpha \in
(0,2\eta\wedge d)$.
\end{example}

Let us now introduce the set of hypotheses to be assumed in order to prove
existence and uniqueness of solution for (\ref{1.4}) and H\"older
properties for its paths. Given $\gamma_1, \gamma_2 >0$, we denote
by $\cc^{\gamma_1, \gamma_2 }([0,t]\times \R^d)$, the space of
real-valued functions defined on $[0,t]\times \R^d$, jointly H\"older
continuous of order $\gamma_1$ in its first variable and $\gamma_2$
in its second one. 

\begin{description}
\item[(H1)] 
For any $i,j=1,\cdots,n$, 
$a^{i,j}, b^i: \Omega\times \R_+\times \R^d \rightarrow \R $ are
$\mathcal{P}\otimes \mathcal{B}(\R^d)$-measurable such that, for any $ \omega\in \Omega$ a.s. and
$t\leq 0$, $a^{i,j}(\omega,t,\cdot)\in \cc^{\alpha}(\R^d)$,
$\alpha\in (\frac{3}{2},2)$, $b^i(\omega,t,\cdot)\in
\cc^{0,1}(\R^d)$, and
\[
\sup_{t\le 0}
\left[ \|a(t,\cdot)\|_{\cc^\alpha}
+ \|b(t,\cdot)\|_{\cc^{0,1}}
\right] \leq k.
\]
There exist $K, \delta>0$,
such that for any $\omega\in \Omega$ a.s., $t\geq 0$,
$x,\lambda\in \R^d$,
\[
\delta |\lambda|^2 \leq a^{i,j}(t,x) \lambda^i \lambda^j \leq K
|\lambda|^2 .
\]

\item[(H2)]
$f,h: \Omega\times\R_+\times \R^d \times \R \rightarrow \R $
are such that, for any $x$ and $u$,
$f(\cdot,x,u)$,
$h(\cdot,x,u)$ are predictable and
\[
\sup_{(\omega,t,x)\in\Omega\times\R_+\R^d}\big[|f(t,x,u)-f(t,x,v)|
+|h(t,x,u)-h(t,x,v)|\big] \leq k |u-v|,
\]
for some positive constant $k$, a.s.
\end{description}

We recall that, for $\alpha\in(\frac{3}{2}, 2)$, $\cc^\alpha$ is the space of continuously 
differentiable functions whose partial derivatives of first order are $\{\alpha\}$- H\"older
continuous, where $\alpha=[\alpha] +\{\alpha\} $, $[\alpha]$ meaning the integer part of
$\alpha$ (see \cite{t1});  
$\cc^{0,1}$ is the space of Lipschitz continuous functions.

In the proof of the next theorem we will use the following Remark 5.5 of \cite{k2}:

\begin{quote}
\textit{For any $u\in H_p^{n+2}$, $ m \in [n,n+2]$ and $\varepsilon > 0$,  we have 
\[
\begin{array}{ll}
\|u\|_{m,p}
&
\leq N \|u\|^{\theta}_{n+2,p}
\|u\|^{1-\theta}_{n,p}
\\ \\
&
\leq
N \theta \varepsilon \|u\|_{n+2,p} +
N (1-\theta) \varepsilon^{-\frac{\theta}{1-\theta}}
\|u\|_{n,p} ,
\end{array}
\]
where $ \theta = \frac{m-n}{2}$ and $N$ depends only on $d$, $n$,
$m$ and $p$.}
\end{quote}

In the following theorem $\tau$ denotes a fixed stopping time with respect to the filtration
$\{\mathcal{F}_t, t\ge 0\}$ defined in Section 2.

\begin{theorem}
\label{th4}
Suppose that there exists
$ \eta \in (\frac{1}{2},1)$
such that
\[
\nu_{\eta,d} = \|R_{\eta,d}\|^2_\ch < + \infty .
\]
We also assume that, for some $p\in [2,+\infty)$
the following conditions are satisfied:
\begin{description}
\item[(a)]
$u_0\in L_p(\Omega, \cf_0, H_p^{1-\eta-\frac{2}{p}}) $,
\item[(b)]
\begin{equation}
\label{1.7}
I^p(\tau)=E\left[\int_0^\tau dt \Big(
\|f(t,\cdot,0)\|_{-1-\eta,p}^p
+
\|\overline{h}(t,\cdot,0)\|_{p}^p \Big) \right]
< +\infty ,
\end{equation}
where
\begin{equation}
\label{1.8}
\overline{h}(t,x,0) :=
\|R_{\eta,d}(x-\cdot) h(t,\cdot,0)\|_\ch .
\end{equation}
\end{description}
Then, in the space $\ch_p^{1-\eta}(\tau)$, equation (\ref{1.4})
with the initial condition $u_0$ and coefficients satisfying
(H1),(H2) posseses a unique solution $u$. Moreover,
\begin{equation}
\label{1.81}
\|u\|_{\ch_p^{1-\eta}(\tau)} \leq
C \left(
I(\tau) +
\Big(E(\|u_0\|_{1-\eta-\frac{2}{p}}^p)^{\frac{1}{p}}
\Big) \right) ,
\end{equation}
where the constant $C$ depends on
$\eta$, $d$, $\alpha$, $p$, $\delta$, $k$ and $\tau$.

In addition, if conditions (a), (b) are satisfied for any $p\geq
2$ then, the trajectories of $u$ 
belong to the space of H\"older continuous functions
$\cc^{\gamma_1,\gamma_2}([0,\tau]\times\R^d)$, a.s.
with $\gamma_1\in(0,\frac{1-\eta}{2})$, $\gamma_2\in(0,1-\eta)$.
\end{theorem}

\begin{proof}

The existence and uniqueness of solution will follow by applying
 Theorem \ref{th1} to
\[
f(t,x,u):=b^i(t,x) u_{x^i}(t,x) + f(t,x,u(t,x)) ,
\]
\[
g^k(t,x,u):=h(t,x,u(t,x)) v_k(x) ,
\]
and by taking $n=-(1+\eta)$.
In fact, we will check that the hypotheses (A1) and (A2) are satisfied.

Since $n\in(-2,-\frac{3}{2})$, we shall consider as space
$B^{|n|+\gamma}$, with $\gamma > 0$ and $|n|+\gamma$ not an
integer, the space $\cc^\alpha(\R^d)$, with $\alpha\in
(\frac{3}{2},2)$.

Set $\overline{f}(t,x,u) = b^i(t,x) u_{x^i}(t,x) + f(t,x,u(t,x))$;
we have to check the following conditions for $n=-(1+\eta)$ (see Assumption {\bf(A2)} before):
\begin{description}
\item[(1)] For any $u\in H_p^{n+2}$, $\{\overline{f}(t,\cdot,u),
t\geq 0\}$ is a predictable process with values on $H^n_p$.
\item[(2)] $\overline{f}(\cdot,\ast,0)\in \HH_p^n$, a.s.
\item[(3)] $\overline{f}$ is a continuous function in $u$ a.s.
\item[(4)] For any $\varepsilon>0$, there exists $K_\varepsilon$
such that, for every $u, v\in H_p^{n+2}, t, \omega$,
\[
\|\overline{f}(t,\cdot,u)-\overline{f}(t,\cdot,v)\|_{n,p}
\leq \varepsilon \|u-v\|_{n+2,p}+K_\varepsilon \|u-v\|_{n,p}.
\]
\end{description}

The predictability of $\overline{f}$ clearly follows from the same
property of $b$ and $f$.

Let $u\in H_p^{n+2}$; then $u_{x^i}\in H_p^{n+1}$.
Notice that,
since $|n+1|\in (\frac{1}{2},1)$, the space $B^{|n+1|+\gamma}$
coincides with the space of the $\alpha$-H\"older continuous
functions for some $\alpha\in (0,1)$.
Since $\cc^{0,1}(\R^d) \subset \cc^\alpha(\R^d)$,
Lemma 5.2 in Krylov \cite{k2} applied to
$b$ and $u$ yields
\begin{equation}
\label{1.9}
\|b^i u_{x^i}\|_{n+1,p} \leq
\|b\|_{B^{|n+1|+\gamma}} \|u_{x^i}\|_{n+1,p}< +\infty .
\end{equation}
Thus $b^i u_{x^i}\in H_p^n$.

Moreover, $u\in L_p$ and the Lipschitz condition of $f$ with
respect to $u$ implies $f(u)-f(0)\in L_p \subset H_p^n$. By
(\ref{1.7}), $f(t,\cdot,0)\in H_p^n$ a.e. on $\ppa0,\tau\ppa$.
Therefore,  {\bf (2)} holds. In addition
\[
|f(t,\cdot,u)|\leq k|u| + |f(t,\cdot,0)|.
\]
This proves $f(t,\cdot,u)\in H_p^n$ and thus
$\overline{f}(t,\cdot,u)\in H_p^n$ as well.

Since $n<-1$, applying again Lemma 5.2 in \cite{k2}, we have
\begin{equation}
\label{1.10}
\begin{array}{ll}
\|b^i u_{x^i}\|_{n,p}
&
\leq
\|b^i u_{x^i}\|_{-1,p} \leq
\|b\|_{\cc^{0,1}} \|u_{x^i}\|_{-1,p}
\\ \\
&
\leq K \|u\|_p = K \|u\|_{n+2+\eta-1,p},
\end{array}
\end{equation}
where in the last identity we have used that $n+2+\eta-1=0$. This
fact, together with the Lipschitz property of $f$ with respect to
$u$, prove {\bf (3)}.

We also have
\begin{equation}
\label{1.11}
\|f(t,\cdot,u)-f(t,\cdot,v)\|_{n,p} \leq
\|f(t,\cdot,u)-f(t,\cdot,v)\|_{p} \leq
K \|u-v\|_p.
\end{equation}
Then
\[
\|\overline{f}(t,\cdot,u)-\overline{f}(t,\cdot,v)\|_{n,p} \leq
A + B,
\]
where
\[
A=\|f(t,\cdot,u)-f(t,\cdot,v)\|_{n,p} \leq
K \|u-v\|_p,
\]
\[
B=\|b^i u_{x^i}-b^i v_{x^i}\|_{n,p} \leq
K \|u-v\|_p,
\]
by (\ref{1.11}) and (\ref{1.10}), respectively.

We now apply the above quoted Remark 5.5 of Krylov
\cite{k2} to $m=n+2+\eta-1=0$.
Notice that $\theta=-\frac{n}{2}>0$,
$1-\theta=\frac{2+n}{2}>0$,
$-\frac{\theta}{1-\theta}=\frac{n}{2+n}<0$.
This yields property {\bf (4)}.

Concerning the coefficient
$g(t,x,u)=\{h(t,x,u(t,x) v_k(x)\}_{k\geq 0}$ we have to check
first of all that, for any $u\in H_p^{n+2}$, $\{g(t,\cdot,u),
t\geq 0\}$ is a predictable process with values on
$H^{n+1}_p(l_2)$. This is a simple consequence of the fact that
$h$ is predictable and $v_k(x)$ is deterministic. Moreover,
since $h$ is Lipschitz, it is immediate also to prove that $g$ is
a.s. continuous in $u$.

\no
Let us now prove that $g(\cdot,\ast,0)\in
\HH_p^{n+1}(\tau,l_2)$.
Since $n+1 = -\eta$, Lemma \ref{lm1} yields
\[
\begin{array}{ll}
\|g(t,\cdot,u)\|_{n+1,p} & \leq
\|g(t,\cdot,u)-g(t,\cdot,0)\|_{n+1,p} + \|g(t,\cdot,0)\|_{n+1,p}
\\ \\
& \leq \|h(t,\cdot,u)-h(t,\cdot,0)\|_{p} + \|h(t,\cdot,0)\|_{p}
\\ \\
& \leq |u| + \|h(t,\cdot,0)\|_{p} < +\infty .
\end{array}
\]
and
\[
\|g(t,\cdot,0)\|_{n+1,p} = \|\overline{h}(t,\cdot,0)\|_{p},
\]
with $\overline{h}(t,\cdot,0)$ defined in (\ref{1.8}).
Thus
\[
\begin{array}{ll}
\|g(t,\cdot,0)\|^p_{\HH_p^{n+1}(\tau,l_2)} & {\ds =
E\Big(\int_0^\tau \|g(t,\ast,0)\|^p_{n+1,p} dt \Big)}
\\ \\
& {\ds = E\Big(\int_0^\tau \|\overline{h}(t,\ast,0)\|^p_{p} dt
\Big) < +\infty ,}
\end{array}
\]
by (\ref{1.7}) and we obtain $g(\cdot,\ast,0)\in
\HH_p^{n+1}(\tau,l_2)$.

It remains to check that for any $\varepsilon>0$, there exists a
constant $K_\varepsilon$ such that, for any $u, v\in H_p^{n+2}, t,
\omega$,
\[
\|g(t,\cdot,u)-g(t,\cdot,v)\|_{n+1,p}
\leq \varepsilon \|u-v\|_{n+2,p}+K_\varepsilon \|u-v\|_{n,p}.
\]
Applying again Lemma \ref{lm1} and the Lipschitz property of $h$,
yield
\[
\begin{array}{ll}
\|g(t,\cdot,u)-g(t,\cdot,v)\|_{n+1,p}
&
\leq
C \|h(t,\cdot,u)-h(t,\cdot,v)\|_{p}
\\ \\
&
\leq
K\|u-v\|_p = K \|u-v\|_{n+2+\eta-1,p}.
\end{array}
\]
Then the above property follows, as for $\overline{f}$, from the
above-mentioned Remark 5.5 in \cite{k2}.

This finishes the
proof of the existence and uniqueness of solution for equation
(\ref{1.4}) in the space $\ch_p^{1-\eta}(\tau)$, and of the bound
(\ref{1.81}).

Let us now check the H\"older continuity of the trajectories of
the solution, in a similar way as in Remark 8.7 in  \cite{k2}. Let $p>2$,
$\frac{1}{2}>\beta>\alpha>\frac{1}{p}$. Then by
Theorem 7.2 in  \cite{k2}, a.s. 
\[
u\in \cc^{\alpha-1/p}([0,\tau], H_p^{1-\eta-2\beta}). 
\]

The space $H_p^{1-\eta-2\beta}$ is
embedded into $\cc^\gamma(\R^d)$ for
$\gamma<1-\eta-2\beta-\frac{d}{p}$, whenever
$1-\eta-2\beta-\frac{d}{p}>0$ (see for instance Theorem E.12 of \cite{shi}).
Thus, taking $p$ big enough and
$\alpha, \beta$ small, we prove that $u$ is $\gamma_2$- H\"older
continuous in $x$, with $\gamma_2<1-\eta$, uniformly in $t$. On
the other hand, the  conditions $\beta< \frac{1}{2}$, $1-\eta-2\beta-\frac{d}{p}>0$
are simultaneously satisfied for  $p$ big enough whenever
$\beta<\frac{1-\eta}{2}\wedge \frac{1}{2}=\frac{1-\eta}{2}$.
Thus $u$ is H\"older continuous in $t$ of order
$\gamma_1<\frac{1-\eta}{2}$, uniformly in x.

This finishes the proof of the theorem.
\end{proof}

\begin{remark}
\label{rem4}
The assumptions of Theorem \ref{th4} ensuring H\"older continuity are satisfied if, for
instance,
$u_0(\omega,\cdot)$ is a.s. a $\cc^\infty$ function with
compact support and
$f(t,x,0)=h(t,x,0)=0$.
\end{remark}


\section{Mild formulation: results on the existence and
uniqueness of a solution}

In this section we consider  the formal expression (\ref{1.4}),
but now, we assume that the coefficients $a, b$ are deterministic.
More precisely, we fix a finite time horizon $T>0$ and 
we assume the following set of assumptions:

\begin{description}
\item[(H1$^\prime$)]
$a^{i,j},b^i: [0,T]\times \R^d \rightarrow \R $,
$i, j=1, \ldots, d$ are $\frac{\alpha}{2}$-H\"older continuous in $t\in[0,T]$,
$\alpha$-H\"older continuous in $x\in\R^d$, for some $\alpha\in(0,1)$.
In addition, for any $\lambda\in \R^d$, there exist $K,\delta>0$ such that
\begin{equation}
\label{2.1} \delta |\lambda|^2 \leq a^{i,j}(t,x) \lambda^i
\lambda^j \leq K |\lambda|^2.
\end{equation}

\item[(H2$^\prime$)]
$f,h:
\Omega\times [0,T]\times \R^d \times \R \rightarrow \R $
are such that, for any $x\in \R^d$ and $u\in \R$,
$f(\cdot,x,u)$,
$h(\cdot,x,u)$ are predictable processes 
satisfying the Lipschitz condition
\begin{align*}
&\sup_{(\omega,t,x)\in\Omega\times [0,T]\times \R^d}
\left[
|f(t,x,u)-f(t,x,v)|
+
|h(t,x,u)-h(t,x,v)|
\right]\\
& \quad \leq k |u-v| ,
\end{align*}
for any $u,v\in \R$.
\end{description}
Following classical approaches on spde's, one can think of equation
(\ref{1.4}) with the initial condition $ u(0,x)=u_0(x) $ as a
stochastic Cauchy problem 
\begin{equation}
\label{2.2}
\left\{
\begin{array}{ll}
\cl u(t,x)= f(t,x,u(t,x))+h(t,x,u(t,x)) F(dt,dx)
\\ \\
u(0,x)=u_0(x)
\end{array}
\right.
\end{equation}
where $\cl$ is the second order operator with coefficients
depending on $t$ and $x$, acting on functions defined on
$[0,T]\times \R^d$, given by
\begin{equation}
\label{op}
\cl=\frac{\partial}{\partial t} -
\sum_{i, j=1}^d a^{i, j}(t,x) \partial^2_{x_i x_j}
- \sum_{i=1}^d b^{i}(t,x) \partial_{x_i} .
\end{equation}
By virtue of (\ref{2.1}) the operator $\cl$ is uniformly parabolic
in $[0,T]\times \R^d$ (\cite{lsu}, page 11)).

Let $G(t,x;s,y)$ be the fundamental solution of $\cl u=0$. $G$ is
a function defined on $[0,T]\times \R^d\times [0,T]\times \R^d
\cap \{(s,t): 0\leq s \leq t \leq T \}$. Under the above
assumptions on the coefficients of $\cl$, $G$ is continuous in all
its variables and for any fixed $s\in [0,T], y \in \R^d$,
$G(\cdot, \ast,;s,y)$ is twice continuously differentiable in $x$,
once continuously differentiable in $t$ and satisfies the
estimates
\begin{equation}
\label{2.3}
\left| \partial_x^\mu \partial_t^\nu G(t,x;s,y)\right|
\leq C (t-s)^{-\frac{d+|\mu|+2\nu}{2}}
\exp\left(-c \frac{|x-y|^2}{t-s}\right)
\end{equation}
where $\mu=(\mu_1,\ldots, \mu_d)\in \N^d, \nu\in \N$, $
|\mu|+2\nu\leq 2$, with $|\mu|=\sum_{j=1}^d \mu_j$ (see (13.3), page 376 in \cite{lsu}). 
Moreover, $G$ is a positive function (Theorem 11 in \cite{fri}).

Let us now introduce the notion of {\it mild solution}. A predictable
stochastic process $\{u^M(t,x), (t,x)\in [0,T]\times \R^d \}$ is
said to be a mild solution to the stochastic Cauchy problem
(\ref{2.2}) if it satisfies the equation
\begin{align}
u^M(t,x)& = \int_{\R^d} dy G(t,x;0,y) u_0(y)\nonumber\\
& + \int_0^t \int_{\R^d} F(ds,dy) G(t,x;s,y) h(s,y,u^M(s,y))\nonumber\\
& + \int_0^t ds \int_{\R^d} dy G(t,x;s,y) f(s,y,u^M(s,y))\label{2.4} .
\end{align}

Notice that, in order to give a rigourous meaning to
equation (\ref{2.4}), we must specify the space where
the solution belongs to.

Using the CONS $\{e_k, k\geq 0\}$ of $\ch$ introduced in Section
2, the stochastic integral in (\ref{2.4}) can also be written as
\[
\sum_{k=0}^{\infty}
\int_0^t  W^k(ds)
\langle G(t,x;s,\cdot) h(s,\cdot,u^M(s,\cdot)),e_k \rangle,
\]
with $ W^k(t)=\int_0^t \int_{\R^d} F(ds,dy) e_k(y)$.

In the sequel, we denote by $G_0(t,x)$ the $d$-dimensional Gaussian
density, zero mean, with variance $t \,\mbox{Id}_n$. Notice that
(\ref{2.3}) implies
\[
|G(t,x;s,y)|=G(t,x;s,y)\leq C_1 G_0(C_2 (t-s),(x-y)),
\]
for some positive constants $C_1, C_2$.

For any $p\in [2,\infty)$, 
let $\cb_p$ be the Banach space of real-valued predictable
processes such that
\[
\sup_{(t,x)\in [0,T]\times \R^d} E(|u(t,x)|^p)< \infty
.
\]

\begin{theorem}
\label{thm2.1}
Fix $p\in[2,\infty)$. Assume (H1$^\prime$),
(H2$^\prime$) and 
\begin{equation}
\label{2.5}
\int_0^T ds \sup_{y\in \R^d} E\left(|h(s,y,0)|^p+|f(s,y,0)|^p\right) < \infty .
\end{equation}
Suppose also that $\|u_0\|_\infty <C$ and moreover,
\begin{equation}
\label{2.6}
\int_{\R^d}
\frac{\mu(d\xi)}{1+|\xi|^2} < \infty.
\end{equation}
Then, there exists a unique stochastic processes $\{u^M(t,x), (t,x)\in[0,T]\times \R^d\}$
belonging to  $\cb_p$ and satisfying  (\ref{2.4}).
\end{theorem}

\begin{proof}
Following classical ideas (see e.g. \cite{sv1}), we show that the map
\begin{align}
\ct u(t,x) & = \int_{\R^d} dy G(t,x;0,y) u_0(y)\nonumber\\
& + \int_0^t \int_{\R^d} F(ds,dy) G(t,x;s,y) h(s,y,u(s,y))\nonumber\\
& + \int_0^t ds \int_{\R^d} dy G(t,x;s,y) f(s,y,u(s,y))\label{2.7}
\end{align}
possesses a unique fixed point in $\cb_p$.

First, we show that (\ref{2.7}) defines a map
$\ct:\cb_p \rightarrow \cb_p$.
Indeed,
\[
E\Big(|\ct u(t,x)|^p\Big) \leq C \Big(T_1(t,x) +
T_2(t,x) + T_3(t,x)\Big),
\]
with
\[
\begin{array}{l}
{\ds T_1(t,x) = \left| \int_{\R^d} dy G(t,x;0,y) u_0(y)
\right|^p,}
\\
{\ds T_2(t,x) = E\left(\left|\int_0^t \int_{\R^d} F(ds,dy) G(t,x;s,y)
h(s,y,u(s,y)) \right|^p\right ),}
\\
{\ds T_3(t,x) = E\left (\left|\int_0^t ds \int_{\R^d} dy G(t,x;s,y)
f(s,y,u(s,y)) \right|^p\right) .}
\end{array}
\]
Clearly,
\[
\sup_{(t,x)\in [0,T]\times \R^d}T_1(t,x) \leq \|u_0\|^p_\infty \sup_{(t,x)\in [0,T]\times \R^d}
\left| \int_{\R^d} G(t,x;0,y) dy \right|^p < \infty.
\]

Notice that

\begin{align}
&\int_{\R^d} \Gamma(dz) \ \big(G_0 (t-s,x-\cdot)
\ast G_0 (t-s,x-\cdot)\big)(z)\nonumber\\
&= \int_0^t ds \int_{\R^d} \mu(d\xi) |\cf G_0(s,\cdot)(\xi)|^2
\leq \int_{\R^d} \frac{\mu(d\xi)}{1+|\xi|^2} <\infty,\label{2.7.1}
\end{align}
which implies
\[
\begin{array}{l}
{\ds \int_0^t ds \int_{\R^d} \Gamma(dz) 
\big(G(t,x;s,\cdot)*\tilde G(t,x;s,\cdot)\big)(z)}
\\
{\ds \leq \int_0^t ds \int_{\R^d} \Gamma(dz) \ \big(G_0(t-s,x-\cdot)
\ast G_0(t-s,x-\cdot)\big)(z)}
\\
{\ds \leq \int_0^t ds \int_{\R^d} \frac{\mu(d\xi)}{1+|\xi|^2} 
< \infty.}
\end{array}
\]
Consequently
\begin{equation}
\label{2.7.2}
\sup_{0\leq t \leq T}
\int_0^t ds \int_{\R^d} \Gamma(dz)
\big(G(t,x;s,\cdot)*\tilde G(t,x;s,\cdot)\big)(z) < \infty .
\end{equation}
Burkholder's and H\"older's inequalities yield
\[
\begin{array}{rl}
T_2(t,x) \leq & {\ds E\Big( \int_0^t ds \int_{\R^d} \Gamma(dz)
\int_{\R^d} dy G(t,x;s,y) h(s,y,u(s,y))}
\\
&
\times G(t,x;s,y-z) h(s,y-z,u(s,y-z))\Big)^{\frac{p}{2}}
\\
\leq & {\ds \left( \int_0^t ds \int_{\R^d} \Gamma(dz) 
\big( G(t,x;s,\cdot)*\tilde G(t,x;s,\cdot)\big)(z)\right)^{\frac{p}{2}-1}}
\\
& {\ds \times \int_0^t ds \int_{\R^d} \Gamma(dz) \int_{\R^d} dy \
G(t,x;s,y) G(t,x;s,y-z)}
\\
& {\ds \times E\Big( |h(s,y-z,u(s,y-z))|^{\frac{p}{2}}
|h(s,y,u(s,y))|^{\frac{p}{2}} \Big).}
\end{array}
\]

Hence, by virtue of (\ref{2.7.1}), (\ref{2.7.2}) and the Lipschitz continuity of $h$,
we obtain
\begin{align*}
&\sup_{(t,x)\in [0,T]\times \R^d}T_2(t,x) \leq  C 
 \sup_{(t,x)\in [0,T]\times \R^d}\int_0^t ds \sup_y E\left(|u(s,y)|^p +
|h(s,y,0)|^p \right)\\
& \quad \times \int_{\R^d} \Gamma(dz) \ \ G_0(t-s,x-\cdot) \ast
G_0(t-s,x-\cdot) < \infty,
\end{align*}
Similarly,
\[
\sup_{(t,x)\in [0,T]\times \R^d} T_3(t,x) < \infty .
\]
Therefore $ \ct u(t,x) $ belongs to the Banach space
$\cb_p$.

To conclude the proof, one can follow the arguments of the proof
of Proposition 3 in  \cite{sv1} (see
also the end of the proof of the next theorem).

\end{proof}

In order to compare mild and weak solutions, we need a 
new version of existence of mild solution.

\begin{theorem}
\label{thm2.2}
Fix $p\in[2,\infty)$. Assume (H1$^\prime$), (H2$^\prime$) and
\begin{equation}
\label{2.8}
E\left[\int_0^T ds
\left(
\|h(s,\cdot,0)\|_p^p
+
\|f(s,y,0)\|_p^p
\right)\right] < \infty .
\end{equation}
Suppose also that $u_0 \in L_p(\R^d)$ and
\[
\int_{\R^d}
\frac{\mu(d\xi)}{1+|\xi|^2} < \infty.
\]
Then, there exists a unique stochastic process
$\{u^M(t,x), (t,x)\in[0,T]\times\R^d\}$ that belongs to $L_p(\Omega\times [0,T]; L_p(\R^d))$
and satisfies (\ref{2.4}).
\end{theorem}

\begin{proof}
We shall divide the proof into two steps, proving first that the
map defined by (\ref{2.7}) on  $L_p(\Omega\times
[0,T]; L_p(\R^d))$, takes values on the same space and secondly, that
this map is a contraction.

{\bf Step 1}: Let $u\in L_p(\Omega\times [0,T]; L_p(\R^d))$; then,
\[
E\left(\int_0^T dt \int_{\R^d} dx |\ct u(t,x)|^p\right)
\leq C(T_1 + T_2 + T_3),
\]
with
\[
\begin{array}{l}
{\ds T_1 = \int_0^T dt \int_{\R^d} dx \left|\int_{\R^d} dy
G(t,x;0,y) u_0(y)\right|^p,}
\\
{\ds T_2 = \int_0^T dt \int_{\R^d} dx E \left(\left|\int_0^t \int_{\R^d}  F(ds,dy)
G(t,x;s,y) h(s,y,u(s,y)) \right|^p\right),}
\\
{\ds T_3 = \int_0^T dt \int_{\R^d} dx E \left(\left|\int_0^t ds
\int_{\R^d} dy G(t,x;s,y) f(s,y,u(s,y)) \right|\right)^p .}
\end{array}
\]
H\"older's inequality, the properties of $G$ and Fubini's
theorem yield
\begin{equation}
\label{2.8.5}
\begin{array}{rl}
T_1
&
{\ds \leq \int_0^T dt \int_{\R^d} dx \left(\int_{\R^d} dy
G(t,x;0,y) \right)^{p-1} \int_{\R^d} dy G(t,x;0,y)
\left|u_0(y)\right|^p}
\\
& {\ds \leq \sup_{(t,x)\in [0,T]\times \R^d} \left(\int_{\R^d} dy
G(t,x;0,y) \right)^{p-1} \int_0^T dt \int_{\R^d} dy
\left|u_0(y)\right|^p}
\\
& {\ds \quad \times \int_{\R^d} dx G(t,x;0,y)}
\\
& {\ds \leq C \sup_{(t,x)\in [0,T]\times \R^d} \left(\int_{\R^d}
dy G_0(t,x-y) \right)^{p} \|u_0\|^p_p < \infty .}
\end{array}
\end{equation}
To deal with $T_2$, we apply Burkholder's inequality, then
H\"older's inequality; we obtain
\begin{align}
T_2 \leq &  C \int_0^T dt \int_{\R^d} dx E \big(\int_0^t ds
\int_{\R^d} \Gamma(dz) \int_{\R^d} dy G(t,x;s,y) h(s,y,u(s,y))\nonumber\\
&  \times G(t,x;s,y-z) h(s,y-z,u(s,y-z))\big)^{\frac{p}{2}}\nonumber\\
\leq &  C \int_0^T dt \int_{\R^d} dx  E \big(\int_0^t ds
\int_{\R^d} \Gamma(dz) \int_{\R^d} dy G_0(t-s,y-x)\nonumber\\
& \times G_0(t-s,y-z-x) |h(s,y,u(s,y))|
|h(s,y-z,u(s,y-z))| \big)^{\frac{p}{2}}\nonumber\\
\leq & C \int_0^T dt \int_{\R^d} dx \big(\int_0^t ds
\int_{\R^d} \Gamma(dz) 
( G_0(t-s,\cdot-x)*\tilde G_0(t-s,\cdot-x))(z)\big)^{\frac{p}{2}-1}\nonumber\\
& \times \int_0^t ds \int_{\R^d} \Gamma(dz) \int_{\R^d} dy
G_0(t-s,y-x) G_0(t-s,y-z-x)\nonumber\\
&\times E \big( |h(s,y,u(s,y))|^{\frac{p}{2}}
|h(s,y-z,u(s,y-z))|^{\frac{p}{2}} \big) .\label{2.9}
\end{align}
Then, owing to (\ref{2.7.1}),
\[
\begin{array}{rl}
T_2 \leq & {\ds C \int_0^T dt \int_{\R^d} dx \int_0^t ds
\int_{\R^d} \Gamma(dz) \int_{\R^d} dy G_0(t-s,y-x) G_0(t-s,y-z-x)}
\\
& {\ds \times E \left( |h(s,y,u(s,y))|^{\frac{p}{2}}
|h(s,y-z,u(s,y-z))|^{\frac{p}{2}} \right) .}
\end{array}
\]
Since the covariance functional is translation invariant, the last
expression is bounded by

\[
\begin{array}{l}
{\ds C \int_0^T dt \int_{\R^d} dx \int_0^t ds \int_{\R^d}
\Gamma(dz) \int_{\R^d} dy G_0(t-s,y) G_0(t-s,y-z)}
\\
{\ds \ \ \times E \left( |h(s,y+x,u(s,y+x))|^{\frac{p}{2}}
|h(s,y-z+x,u(s,y-z+x))|^{\frac{p}{2}} \right) .}
\end{array}
\]

Applying Fubini's theorem and Schwarz inequality yields
\begin{align}
T_2 &\leq  C \int_0^T dt \int_{\R^d} dx \int_0^t ds
\int_{\R^d} \Gamma(dz) \int_{\R^d} dy G_0(t-s,y) G_0(t-s,y-z)\nonumber\\
&\quad\times E \big( |h(s,y+x,u(s,y+x))|^{\frac{p}{2}}
|h(s,y-z+x,u(s,y-z+x))|^{\frac{p}{2}} \big)\nonumber\\
&\leq  C \int_0^T dt \int_0^t ds \int_{\R^d} \Gamma(dz)
\int_{\R^d} dy G_0(t-s,y) G_0(t-s,y-z)\nonumber\\
& \quad\times E \big( \int_{\R^d} dx
|h(s,y+x,u(s,y+x))|^p\big)^{\frac{1}{2}}\nonumber\\
& \quad\times E \big( \int_{\R^d} dx
|h(s,y-z+x,u(s,y-z+x))|^p\big)^{\frac{1}{2}}\nonumber\\
&= C \int_0^T dt \int_0^t ds  E \big( \int_{\R^d} dx
|h(s,x,u(s,x))|^p\big)\nonumber\\
& \quad \int_{\R^d} \Gamma(dz) (G_0(t-s,\cdot)*\tilde G_0(t-s,\cdot))(z)\label{2.11}
\end{align}

where the last identity holds by the translation invariance
of Lebesgue measure.

The Lipschitz continuity of $h$ yields
\begin{equation}
\label{2.11.5}
E \left(\int_{\R^d} dx |h(s,x,u(s,x)|^p\right)
\leq C
E \left(\|u(s,\cdot)\|^p_p+\|h(s,\cdot,0)\|^p_p\right)
\end{equation}
Hence, by (\ref{2.7.1}) and (\ref{2.8}),
\begin{equation}
\label{2.12}
\begin{array}{rl}
T_2 \leq & {\ds C_1 \int_0^T dt \int_0^t ds E \left(
\|u(s,\cdot)\|^p_p \right) + C_2 \int_0^T dt \int_0^t ds
E \left( \|h(s,\cdot,0)\|^p_p \right)}
\\
\leq & {\ds C_1 \|u\|_{L_p(\Omega\times
[0,T];L_p(\R^d))} + C_3 .}
\end{array}
\end{equation}
The analysis of $T_3$ is simpler.
Indeed, H\"older's inequality
implies
\begin{align*}
T_3 &\leq \int_0^T dt \int_{\R^d} dx \left( \int_0^t ds
 \int_{\R^d} dy\, G(t,x;s,y)\right)^{p-1}\\
&\quad \times \int_0^t ds \int_{\R^d} dy \,G(t,x;s,y)
E \left(|f(s,y,u(s,y))|^p\right)\\
&\leq  C \int_0^T dt \int_{\R^d} dx \int_0^t ds \int_{\R^d}
dy G(t,x;s,y) E \left(|f(s,y,u(s,y))|^p\right),
\end{align*}
since
\begin{equation}
\label{2.12.5} \sup_{0\leq s \leq t \leq T, x\in \R^d} \int_{\R^d}
dy\, G(t,x;s,y) < \infty .
\end{equation}
Then, Fubini's theorem yields
\[
\begin{array}{rl}
T_3 \leq & {\ds C \int_0^T dt \int_0^t ds \int_{\R^d} dy E
\left(|f(s,y,u(s,y))|^p\right) \ \int_{\R^d} dx\, G(t,x;s,y)}
\\
\leq & {\ds C \int_0^T dt \int_0^t ds \int_{\R^d} dy E
\left(|f(s,y,u(s,y))|^p\right) .}
\end{array}
\]
The estimate (\ref{2.11.5}) with $h$ replaced
by $f$ and (\ref{2.8})
imply
\begin{equation}
\label{2.13}
\begin{array}{rl}
T_3 \leq & {\ds C \int_0^T dt \int_0^t ds \left\{ E \left(
\|u(s,\cdot)\|^p_p \right) +
\|f(s,\cdot,0)\|^p_p \right\}}
\\
\leq & {\ds C_1  \|u\|_{L_p(\Omega\times
[0,T];L_p(\R^d))} + C_4 .}
\end{array}
\end{equation}
Then, (\ref{2.8.5}), (\ref{2.12}) and (\ref{2.13}) give
\[
\|\ct u\|_{L_p(\Omega\times [0,T];L_p(\R^d))}
\leq
C_1 \|u\|_{L_p(\Omega\times [0,T];L_p(\R^d))} + C_2.
\]
This completes the proof of Step 1.

{\bf Step 2}:
The mapping $\ct$ has a unique fixed point in
$L_p(\Omega\times [0,T];L_p(\R^d))$.

Indeed, let $u_1, u_2 \in L_p(\Omega\times [0,T];L_p(\R^d))$.
Proceeding as in Step 1 and by virtue of
the Lipschitz property of $f$ and $h$, we
obtain
\[
\|\ct u_1 - \ct u_2\|^p_{L_p(\Omega\times [0,t];L_p(\R^d))}
\leq
C_1
\int_0^t ds \|u_1-u_2\|^p_{L_p(\Omega\times [0,s];L_p(\R^d))},
\]
for any
$0 \leq t \leq T$.

Consequently, for $N$ big enough, the N-th iterate
of $\ct$ is a contraction on $L_p(\Omega\times [0,T];L_p(\R^d))$.

\end{proof}


\section{Equivalence between  weak and mild formulations}

We devote this section to study 
the relationship between the notions of solution
introduced previously, for some particular classes of spde's.
As a consequence, we deduce path properties of the mild solution.
We start by giving an equivalent weak formulation. Then, we compare
the weak and mild formulation when the differential operator is self-adjoint 
and has nonrandom coefficients.

Let us consider equation (\ref{1.4}) written in terms of the sequence $\{W^k, k\ge 0\}$ of
independent Brownian motions, that is
\begin{align}
du(t,x)& =\big[a^{i,j}(t,x) u_{x^i,x^j} (t,x) + b^i(t,x) u_{x^i} (t,x)\big]\nonumber\\
& + f(t,x,u(t,x)) dt + g^k(t,x,u(t,x)) W^k(dt),\label{3.0.0}
\end{align}
$t\in[0,T]$, with initial condition $u(0,\cdot)= u_0$. 

We have proved in Theorem \ref{th4} the existence of a unique function-valued stochastic
process $\{u(t), t\in[0,T]\}$ satisfying 
\begin{align}
&(u(t,\cdot),\phi) = (u_0,\phi) + \int_0^t ds
\big(a^{i,j}(s,\cdot)u_{x^i,x^j}(s,\cdot) +(b^i(s,\cdot)u_{x^i}(s,\cdot),\phi\big)\nonumber\\
&+ \int_0^t ds (f(s,\cdot,u(s,\cdot)),\phi) + \int_0^t
W^k(ds) (g^k(s,\cdot,u(s,\cdot)),\phi),\label{3.0}
\end{align}
for all $\phi \in \cc_0^\infty(\R^d)$,
with the pairing $(\cdot,\cdot)$ given in (\ref{1.2}).
We shall say that the process $u$ is a {\it weak solution} of equation
(\ref{3.0}).

The next proposition establishes the equivalence between testing
against functions depending on $x$ and functions depending on $t$ and
$x$. To fix the
notation, denote by $\cc^{1,2}_{t,x;0}$ the space of functions
$f:[0,T]\times \R^d \rightarrow \R$ of class $\cc^1$ in $t$, 
$\cc^2$ in $x$, with compact support. 

\begin{proposition}
\label{prop3.1} We assume that the assumptions of Theorem \ref{th4} are satisfied.
The  stochastic process
$u$ is a  weak solution
if and only if for any function 
$\Phi\in\cc^{1,2}_{t,x;0}$, the following identity holds:
\begin{align}
& (u(t,\cdot),\Phi(t,\cdot)) = (u_0,\Phi(0,\cdot)) + \int_0^t
ds (u(s,\cdot),\partial_s \Phi(s,\cdot))\nonumber\\
&\quad + \int_0^t ds (a^{i,j}(s,\cdot) u_{x^i,x^j}(s,\cdot) + b^i(s,\cdot) u_{x^i}(s,\cdot), \Phi(s,\cdot))\nonumber\\
&\quad + \int_0^t ds (f(s,\cdot,u(s,\cdot)),\Phi(s,\cdot)) +
\int_0^t W^k(ds) (g^k(s,\cdot,u(s,\cdot)),\Phi(s,\cdot)) .\label{3.1}
\end{align}
\end{proposition}

\begin{proof}
The ``only if'' part is trivial. To complete the proof, we proceed into two steps.

{\bf Step 1}:  Let us prove the result
in the case where
\[
\Phi(t,x) = \varphi(t) \phi(x),
\]
with $\varphi \in \cc^1([0,T])$
and $\phi \in \cc^2_0(\R^d)$.

In equation (\ref{3.0}), we set
$t=\sigma$, multiply each term by $\varphi^\prime(\sigma)$ and
then integrate on $(0,t)$ with respect to $\sigma$.
We obtain

\begin{align*}
&\int_0^t d\sigma \varphi^\prime(\sigma) \big(u(\sigma,\cdot),\phi\big)\\
& = \varphi(t) \big( u_0,\phi\big) -
\varphi(0)\big( u_0,\phi\big)\\
& + \int_0^t d\sigma \varphi^\prime(\sigma) \int_0^\sigma
ds \big(a^{i,j}(s,\cdot) u_{x^i,x^j}(s,\cdot) + b^i(s,\cdot)u_{x^i}(s,\cdot), \phi\big)\\
& + \int_0^t d\sigma \varphi^\prime(\sigma) \int_0^\sigma
ds\big(f(s,\cdot,u(s,\cdot)), \phi\big)\\
&+ \int_0^t d\sigma \varphi^\prime(\sigma) \int_0^\sigma
 W^k(ds) \big( g^k(s,\cdot,u(s,\cdot)), \phi\big)
\end{align*}
Set
\begin{align*}
I_1& =\int_0^t d\sigma \varphi^\prime(\sigma) \int_0^\sigma
ds \\
&\quad\times \big(a^{i,j}(s,\cdot) u_{x^i,x^j}(s,\cdot) + b^i(s,\cdot)u_{x^i}(s,\cdot), \phi\big),\\
I_2& = \int_0^t d\sigma \varphi^\prime(\sigma) \int_0^\sigma
ds\big(f(s,\cdot,u(s,\cdot)), \phi\big), \\
I_3& = \int_0^t d\sigma \varphi^\prime(\sigma) \int_0^\sigma
 W^k(ds) \big( g^k(s,\cdot,u(s,\cdot)), \phi\big),
\end{align*}

Integrating by parts we obtain,
\begin{align*}
I_1& = I_1^\prime - \int_0^t ds \varphi(s) \big(a^{i,j}(s,\cdot)u_{x^i,x^j}(s,\cdot)
+b^i(s,\cdot)u_{x^i}(s,\cdot),\phi\big),\\
I_2& = I_2^\prime - \int_0^t ds\varphi(s) 
\big(f(s,\cdot,u(s,\cdot)), \phi\big),\\
I_3 & =  I_3^\prime- \int_0^t W^k(ds) \varphi(s) 
\big(g^k(s,\cdot,u(s,\cdot)), \phi\big).
\end{align*}
with

\begin{align*}
 I_1^\prime& = \varphi(t) \int_0^t ds 
\big(a^{i,j}(s,\cdot) u_{x^i,x^j}(s,\cdot)+b^i(s,\cdot) u_{x^i}(s,\cdot),\phi\big),\\
I_2^\prime& = \varphi(t) \int_0^t ds \big(f(s,\cdot,u(s,\cdot)),\phi\big),\\
I_3^\prime& = \varphi(t) \int_0^t W^k(ds) \big(g^k(s,\cdot,u(s,\cdot)), \phi\big).
\end{align*}

Thus, noticing that
\[
\varphi(t)
[(u_0,\phi) +
I^\prime_1 + I^\prime_2 + I^\prime_3]=
\varphi(t)
(u(t,\cdot),\phi),
\]
we obtain
\begin{align*}
&\int_0^t d\sigma\big(u(\sigma,\cdot),\partial_{\sigma}\Phi(\sigma,\cdot)\big) = 
 \int_0^t d\sigma \varphi^\prime(\sigma) \big(u(\sigma,\cdot), \phi\big)\\
& = \varphi(t) (u(t,\cdot),\phi) - \varphi(0)(u_0,\phi)\\
& -\Big\{ \int_0^t d\sigma \varphi(\sigma)
\big(a^{i,j}(\sigma,\cdot)u_{x^i,x^j}(\sigma,\cdot) + b^i(\sigma,\cdot)u_{x^i}(\sigma,\cdot),\phi\big)\\
 &+ \int_0^t d\sigma \varphi(\sigma) \big(f(\sigma,\cdot,u(\sigma,\cdot)), \phi\big)\\
& + \int_0^t W^k(d\sigma) \varphi(\sigma)
\big( g^k(\sigma,\cdot,u(\sigma\cdot)), \phi\big)\Big\}.
\end{align*}
yielding (\ref{3.1}).

{\bf Step 2}:
The function $\Phi$
can be approximated by polynomials
\[
p(t,x) = \sum_{\alpha, \beta \geq0}
c_{\alpha,\beta} \ x^\alpha t^\beta ,
\]
 $\alpha=(\alpha_1, \ldots, \alpha_d)$, in the norm
\[
\|\Phi\| = \sup_{t,x}
\left(
|\Phi(t,x)|+|\partial_t \Phi(t,x)|
+ \sum_{|k|\le 2} |\partial_x^{|k|} \Phi(t,x)|
\right) .
\]
(see e.g. Kirillov and Gvishiani \cite{kg}, pag.77).
Let $[0,T]\times K= \mbox{supp}\, \Phi$ and $p_n$ be a polynomial
defined for $(t,x)\in [0,T]\times K$ such that
$\|\Phi-p_n\| < \frac{1}{n}$.

We have proved in Step 1 that (\ref{3.1}) holds
with $\Phi:=p_n$. Set
$\psi_n(t,x) = p_n(t,x) - \Phi(t,x)$.
We now check that the $L^1(\Omega)$ norm of any term in
(\ref{3.1}), when $\Phi$ is replaced by  $\psi_n(t,x)$, tends to 0 as $n$
tends to infinity. This shall finish the proof of the proposition.

Indeed, let us first prove that 
\begin{equation}
\label{3.2}
\lim_{n\to\infty}E\Big(\big| \big(u_0,\psi_n(0,\cdot)\big)\big|\Big) = 0.
\end{equation}
Assume first $n_0:=1-\eta-\frac{2}{p} \ge 0$. Since,
\begin{equation*}
||\psi_n(0,\cdot)||_{-n_0,q}\le ||\psi_n(0,\cdot)||_q\le ||\psi_n(0,\cdot)|| <\frac{1}{n},
\end{equation*}
for any $q\in(1,\infty)$, H\"older's inequality with $\frac{1}{p}+\frac{1}{q}=1$ yields
\begin{align*}
E\Big(\big | \big(u_0,\psi_n(0,\cdot)\big)\big|\Big)&\le 
E\Big(\int_{\R^d} dx \big |(1-\Delta)^{\frac{n_0}{2}} u_0(x)\big | \big |
(1-\Delta)^{-\frac{n_0}{2}}\psi_n(0,x)\big|\Big)\\
&\le E\big(||u_0||_{n_0,p}\big)\,||\psi_n(0,\cdot)||_{-n_0,q}\\
&\le C \frac{1}{n},
\end{align*}
yielding (\ref{3.2}).

Assume that  $n_0:=1-\eta-\frac{2}{p} < 0$. The restrictions on $\eta$ and
$p$ yield in this case $-n_0\in (0,1)$. Applying the next Lemma \ref{l2} to
$\psi := \psi_n(0,\cdot)$, we obtain
\begin{equation*}
||\psi_n(0,\cdot)||_{-n_0,q} \le C  ||\psi_n||_q < C \frac{1}{n}.
\end{equation*}
Then we can proceed exactly as before and obtain (\ref{3.2}).

Set $m_0= 1-\eta$. Notice that $m_0>0$. Then, as we did before
for the case $n_0\ge 0$,
\begin{align}
E\Big(\big| \int_0^t ds \big(u(s,\cdot), \partial_s\psi_n(s,\cdot)\big)\big| \Big)
&\le \int_0^t ds E\big(||u(s,\cdot)||_{m_0,p}\big)  ||\psi_n(s,\cdot)||_{-m_0,q}\nonumber\\
&\le C \frac{1}{n}.
\label{3.3}
\end{align}
As in the proof of Theorem \ref{th4}, we set $n= -(1+\eta)$. From Lemma 5.2  \cite{k2},
it follows that
\begin{align*}
&E\Big(\int_0^t ds \big(||a^{i,j}(s,\cdot) u_{x^i,x^j}(s,\cdot)||_{n,p}^p
+ ||b^i(s,\cdot) u_{x^i}(s,\cdot)||_{n+1,p}^p\big)\Big)\\
&< \infty.
\end{align*}
Then, since $n+1<0$, by Lemma \ref{l2} we obtain
\begin{align}
&E\Big(\big| \int_0^t ds \big( a^{i,j}(s,\cdot) u_{x^i,x^j}(s,\cdot)+
b^i(s,\cdot) u_{x^i}(s,\cdot), \psi_n(s,\cdot)\big)\big|\Big)\nonumber\\
\le C\frac{1}{n}.\label{3.4}
\end{align}
Following (\ref{1.11}),
\begin{equation*}
||f(s,\cdot,u(s,\cdot))||_{n,p}\le C\big(||u||_p + ||f(s,\cdot,0)||_{n,p}\big).
\end{equation*}
Consequently, by virtue of (\ref{1.7}) and Lemma \ref{l2}
\begin{align}
&E\Big(\big| \int_0^t ds \big(f(s,\cdot,u(s,\cdot)), \psi_n(s,\cdot)\big|\big)\Big)\nonumber\\
&\le \int_0^t ds E\big(||f(s,\cdot,u(s,\cdot))||_{n,p}\big) ||\psi_n(s,\cdot)||_{-n,q}\nonumber\\
&\le C\frac{1}{n}.\label{3.5}
\end{align}
We now deal with the stochastic integral by considering the $L^2(\Omega)$-norm. The
isometry property of the stochastic integral yields,
\begin{align*}
& E\left| \sum_{k=1}^\infty \int_0^t dW_s^k
(g^k(s,\cdot,u(s,\cdot)),\psi_n(s,\cdot)) \right|^2\\
& = \sum_{k=1}^\infty E \int_0^t ds
\left(g^k(s,\cdot,u(s,\cdot)),\psi_n(s,\cdot) \right)^2
\end{align*}

Using Schwarz's and H\"older's inequalities, this last expression is bounded by
\begin{align*}
&\sup_{s\in[0,T]}||(1-\Delta)^{\frac{\eta}{2}}\psi_n(s,\cdot)||_1 ||(1-\Delta)^{-\frac{\eta}{2}}\psi_n(s,\cdot)||_q\\
&\quad\times \int_0^t ds ||(\sum_{k=1}^\infty
|(1-\Delta)^{-\frac{\eta}{2}} g^k(s,\cdot,u(s,\cdot))|^2)^{\frac{1}{2}}||_p^2 ,\\
\end{align*}
(see \cite{k2}, pg. 191).

Hence, Lemma \ref{l2} and the properties of $g$ imply
\begin{align}
& E\left| \sum_{k=1}^\infty \int_0^t dW_s^k
(g^k(s,\cdot,u(s,\cdot)),\psi_n(s,\cdot)) \right|^2\nonumber\\
&\le C \frac{1}{n}\label{3.6}.
\end{align}
With the estimates (\ref{3.2})-(\ref{3.6}), we finish the proof of the proposition.
\end{proof}
\medskip

We conclude the first part of the section giving the auxiliary result that
has been used in the proof of the preceding proposition.

\begin{lemma}
\label{l2}
For any $\psi\in \mathcal{C}_0^\infty(\R^d)$, $\alpha\in(0,1)$ and $q\in (1,\infty)$,
\begin{equation}
\label{3.7}
||(1-\Delta)^\alpha \psi||_q \le C\big( ||\Delta \psi||_q+||\psi||_q\big).
\end{equation}
\end{lemma}
\begin{proof}
For any $t\in(0,\infty)$, $x\in\R^d$, set
 $\varphi_t(x)= \frac{1}{(4\pi t)^{\frac{d}{2}}}\exp\big(-\frac{|x|^2}{4t}\big)$,
$T_t \psi = \psi*\varphi_t$. It is well-known (see for instance \cite{k3}) that
\begin{equation*}
(1-\Delta)^{\alpha}\psi
 = C(\alpha) \int_0^\infty dt\, \frac{e^{-t}T_t\psi - \psi}{t^{\alpha+1}},
\end{equation*}
where $C(\alpha)$ is a positive constant depending only on $\alpha$.

Fix $t_0 > 1$ and set
\begin{equation*}
A_1 = \big\Vert\int_{t_0}^\infty dt \,\frac{e^{-t}T_t\psi - \psi}{t^{\alpha+1}}\big\Vert_q^q.
\end{equation*}
H\"older's inequality yields
\begin{align*}
A_1 &\le \Big(\int_{t_0}^\infty dt\, t^{-\alpha - 1}\Big)^{q-1} \int_{\R^d} dx \\
&\quad\times\int_{t_0}^\infty dt\, t^{-\alpha - 1} \Big(|e^{-t}T_t \psi(x)|^q + |\psi(x)|^q\Big).
\end{align*}
By Young's inequality,
\begin{equation*}
\sup_{t\in[0,T]}||T_t\psi||_q \le ||\psi||_q\sup_{t\in[0,T]}||\varphi_t||_1 = C_0 ||\psi||_q,
\end{equation*}
with $C_0 ={\ds \sup_{t\in[0,T]}||\varphi_t||_1^q}$.
Thus, 
\begin{equation*}
A_1 \le C ||\psi||_q^q,
\end{equation*}
for some positive constant $C$.

Let us now consider the term
\begin{equation*}
A_2 = \big\Vert\int_0^{t_0} dt \frac{e^{-t}T_t\psi - \psi}{t^{\alpha+1}}\big\Vert_q^q.
\end{equation*}
Notice that, for any $x\in\R^d$, the mapping $t\mapsto e^{-t} T_t \psi (x)$ is differentiable
in $(0,t_0)$
and \begin{align*}
\frac{d}{dt} \big(e^{-t} T_t \psi (x)\big) &= e^{-t}\big(\psi * \Delta \varphi_t - \psi * \varphi_t\big)\\
&=e^{-t}\big(\Delta\psi* \varphi_t- \psi*\varphi_t\big).
\end{align*}
Consequently, Young's inequality yields
\begin{equation*}
\sup_{t\in(0,T)} \Big\Vert\frac{d}{dt} \big(e^{-t} T_t \psi \big)\Big\Vert_q
\le \sup_{t\in[0,T]}||\varphi_t||_1 \big( ||\psi||_q + ||\Delta \psi||_q\big).
\end{equation*}
Therefore, H\"older's inequality implies
\begin{align*}
A_2& \le \int_{\R^d} dx \big(\int_0^{t_0} \frac{dt}{t^\alpha}
 \sup_{t\in(0,T)}\big|\frac{d}{dt} \big(e^{-t} T_t \psi (x)\big)\big|\big)^q\\
&\le \int_{\R^d} dx \big (\int_0^{t_0} \frac{dt}{t^\alpha}
 \sup_{t\in(0,T)}\big|\frac{d}{dt} \big(e^{-t} T_t \psi (x)\big|^q\big)\\
&\le  \sup_{t\in[0,T]}||\varphi_t||_1^q \big( ||\psi||_q^q + ||\Delta \psi||_q^q\big)\\
&\le C\big( ||\psi||_q^q + ||\Delta \psi||_q^q\big).
\end{align*}
Since
\begin{equation*}
||(1-\Delta)^\alpha \psi||_q^q \le C(A_1+A_2),
\end{equation*}
the upper bounds obtained before for $A_1$, $A_2$ give (\ref{3.7}) 
\end{proof}
\medskip

We want now to prove that, if $\{u^W(t,x), (t,x)\in [0,T]\times \R^d\}$
satisfies the weak formulation in the sense of (\ref{3.1}) then, it
also satisfies the mild formulation. To obtain this result, we
restrict the class of operators and assume  that $\cl$ given in (\ref{op})
is self-adjoint.

Let us start by proving the following auxiliary result.

\begin{lemma}
\label{lemma4.1}
Fix $v\in \cc^\infty_0(\R^d)$, $t\in [0,T]$.
Assume that $\mbox{\rm supp } v \subset K$, for some compact set $K\subset\R^d$.
Fix $\varepsilon>0$ and define $v^t\in\cl_{t,x,0}^{1,2}$ as follows:
\begin{equation}
\label{4.1}
v^t(s,x) =
\left\{
\begin{array}{ll}
v(x),
&
\mbox{if } s=t,\, x\in \R^d
\\
{\ds \int_{\R^d} dy\, v(y) G(t,x;s,y),} & \mbox{if } s<t,\, x\in
\stackrel{\circ}{K}
\\
0, & \mbox{if } s\ge t,\, x\in K(\varepsilon),
\end{array}
\right.
\end{equation}
where $K(\varepsilon) = \{x\in\R^d: d(x,K) > \varepsilon\}$.

Then, with the assumptions of Theorem \ref{th4}, the following identity holds
\begin{equation}
\label{4.2}
\begin{array}{rl}
(u^W(t,\cdot),v) = & {\ds (u_0,v^t(0,\cdot)) + \int_0^t
ds (f(s,\cdot,u^W(s,\cdot)),v^t(s,\cdot))}
\\
& {\ds + \int_0^t W^k(ds)
(g^k(s,\cdot,u^W(s,\cdot)),v^t(s,\cdot)).}
\end{array}
\end{equation}

\end{lemma}

\begin{proof}
We recall that
$G(t,x;s,y)$  is of class $\cc^{1}$ in  the variable $t$ and 
$\cc^{2}$ in the variable $x$.
This implies
$ v^t \in \cc^{1,2}_{t,x;0}$.
Hence, $v^t$ can be taken as test function in (\ref{3.1})
and consequently (\ref{4.2}) is equivalent to
\begin{align}
\label{4.3}
0& =  \int_0^t ds (u^W(s,\cdot),\partial_s v^t(s,\cdot))\nonumber\\
& + \int_0^t ds (a^{i,j}(s,\cdot) u^W_{x^i,x^j}(s,\cdot) + b^i(s,\cdot) u^W_{x^i}(s,\cdot), v^t(s,\cdot))\nonumber\\
& =  \int_0^t ds (u^W(s,\cdot),\partial_s v^t(s,\cdot))\nonumber\\
& + \int_0^t ds  (u^W(s,\cdot), \partial_{xî,x^j}^2\big(a^{i,j}(s,\cdot) v^t(s,\cdot)\big) -
 \partial_{x^i} \big(b^i(s,\cdot)v^t(s,\cdot)\big).
\end{align}

By the definition of $v^t$ and the properties
of the fundamental solution \break
$G(\cdot,\ast;s,y)$, this is equivalent to
\begin{equation*}
0 = \int_{\R^d} dy\, v(y) \int_0^t ds\,\big(u^W(s,\cdot), 
\cl^\ast_{s,\cdot}G(t,\cdot;s,y)\big],
\end{equation*}
where $ \cl^\ast$ is the adjoint operator of 
$\cl$ (see Friedman
\cite{fri}, pag. 26), that is
\begin{equation*}
\cl^\ast_{t,x} u(t,x) = -\frac{\partial}{\partial t} - \sum_{i,j=1}^d \partial^2_{x^i,x^j}
\big(a^{i,j}(t,x) u(t,x)\big) + \sum_{i=1}^d \partial_{x_i}\big(b^i(y,x)u(t,x)\big).
\end{equation*}
Since $\cl$ is a self-adjoint operator, $G$ is symmetric in $(x,y)$. Consequently
$\cl^\ast_{s,\cdot}G(t,\cdot;s,y) = \cl^\ast_{s,\cdot}G(t,y; s,\cdot)$.

By Theorem 15 in Friedman \cite{fri}, pag. 28, for every fixed $t\in[0,T]$, $y\in\R^d$,
we have
\[
\cl_{s,x}^\ast G(t,y;s,x) = 0 .
\]
This finishes the proof of the lemma.
\end{proof}

\begin{proposition}
\label{proposition4.2} With the hypotheses of Theorem
{\ref{th4}}, let $\{u^W(t,x), (t,x)\in [0,T]\times \R^d\}$ be a
weak solution in the sense of Proposition \ref{prop3.1}. Then, 
for any $x$-a.e.
\begin{align}
\label{4.4}
u^W(t,x)& = \big(u_0, G(t,x;0,\cdot)\big) +\int_0^t ds\, \big(f(s,\cdot,u^W(s,\cdot)), G(t,x;s,\cdot)\big)\nonumber\\
& + \int_0^t W^k(ds) \big( g^k(s,\cdot,u^W(s,\cdot)), G(t,x;s,\cdot)\big).
\end{align}
\end{proposition}

\begin{proof}
For a fixed $t\in [0,T]$, we write the
expression (\ref{3.1}) with
$\Phi(s,x)=v^t(s,x)$, defined in
(\ref{4.1}).
By virtue of Lemma \ref{lemma4.1}
we obtain
\begin{align*}
(u^W(t,\cdot),v) = & \big( u_0,  \int_{\R^d} dy v(y) G(t,\cdot;0,y)\big)\\
& + \int_0^t ds \big( f(s,\cdot,u^W(s,\cdot)), \int_{\R^d} dy v(y) G(t,\cdot;s,y)\big)\\
& + \int_0^t W^k(ds)\big( g^k(s,\cdot,u^W(s,\cdot)), \int_{\R^d} dy v(y) G(t,\cdot;s,y) \big).
\end{align*}

Fubini's theorem implies
\begin{align*}
(u^W(t,\cdot),v) = &  \int_{\R^d} dy v(y) \big[ \big(u_0, G(t\cdot;0,y)\big)\\
&+ \int_0^t ds \big( f(s,\cdot ,u^W(s,\cdot)), G(t,\cdot;s,y)\big)\\
& +  \int_0^t W^k(ds) \big( g^k(s,\cdot ,u^W(s,\cdot)), G(t,\cdot;s,y)\big)\big]
\end{align*}

Consequently, for any $x\in K$, a.e. with respect to Lebesgue measure, 
\begin{align*}
u^W(t,x)& = \big(u_0, G(t,\cdot;0,x)\big) +\int_0^t ds\, \big(f(s,\cdot,u^W(s,\cdot)), G(t,\cdot;s,x)\big)\nonumber\\
& + \int_0^t W^k(ds) \big( g^k(s,\cdot,u^W(s,\cdot)), G(t,\cdot;s,x)\big).
\end{align*}

Since $G(t,x;s,y)=G(t,y;s,x)$,
for any $s\leq t$ and $x,y\in \R^d$, and $K$ is arbitrary,
this is equivalent to (\ref{4.4}). 
Thus, the proof is complete
\end{proof}

The next result, which is the main conclusion of this section, states that, 
if there exists a function-valued solution in the weak sense, then it must coincide
with the mild solution. We need simultaneously the validity of the assumptions of 
Theorems \ref{th4} and \ref{thm2.2}. More precisely, we have the following theorem.

\begin{theorem}
\label{thm4.3}
Suppose that
\begin{enumerate}
\item
The  operator $\mathcal{L}$ defined in (\ref{op}) is self-adjoint and its coefficients are
 deterministic.

\item The functions $a^{i,j}, b^i: [0,T]\times \R^d\to \R$, $i,j = 1,\cdots, d$
are $\frac{\alpha}{2}$-H\"older continuous in $t\in[0,T]$, for some $\alpha\in(0,1)$ and,
for any $t\in[0,T]$, $a^{i,j}(t,\dot)$ belongs to the space $\mathcal{C}^\beta(\R^d)$, for some
$\beta\in(\frac{3}{2}, 2)$, and $b^i(t,\dot)\in\mathcal{C}^{0,1}(\R^d)$.

\item
For any $\lambda\in \R^d$, there exist $K,\delta>0$ such that
\begin{equation*}
 \delta |\lambda|^2 \leq a^{i,j}(t,x) \lambda^i\lambda^j \leq K |\lambda|^2.
\end{equation*}

\item
The coefficients of Equation (\ref{1.4}) are predictable processes
\begin{equation*}
f,h: \Omega\times [0,T]\times \R^d\times\R \rightarrow \R.
\end{equation*}
For some fixed $p\in [2,\infty)$, the following conditions hold:
\begin{description}
\item{(a)}  For any $u,v \in \R$,
\begin{align*}
&\sup_{(\omega,t,x)\in \Omega\times [0,T]\times \R^d}
\left\{
|f(t,x,u)-f(t,x,v)|+|h(t,x,u)-h(t,x,v)|\right\}\\
&\leq k |u-v|.
\end{align*}
\item{(b)}
$E\big(\int_0^T ds \|h(s,\cdot,0)\|^p_p + \|f(s,\cdot,0)\|^p_p\big) < \infty $.
\end{description}
\item
There exists $\eta \in (\frac{1}{2},1)$ such that
\[
\int_{\R^d}
\frac{\mu(d\xi)}{(1+|\xi|^2)^\eta} < \infty .
\]
\item
$u_0\in L_p(\R^d)$.
\end{enumerate}

Then $u^W=u^M$ as processes in
$L_p(\Omega\times [0,T];L_p(\R^d))$.
Consequently,
$\omega$-a.s., $u^W(t,x)=u^M(t,x)$,
a.e. with respect to Lebesgue measure
on $[0,T]\times \R^d$.
\end{theorem}

\begin{remark}
\begin{enumerate}
\item
The above hypothesis 5 implies that $||R_{\eta,d}||_{\mathcal{H}} < \infty$.

Under this condition, we
have proved in Lemma \ref{lm1}, that
$\|\overline{h}\|_p \leq C \|h\|_p $,
where
$\overline{h}(x)=\|R_{\eta,d}(x-\cdot)h\|_{\ch}$.
Therefore the assumption 4 (b) implies
\[
E\int_0^T ds \|\overline{h}(s,\cdot,0)\|_p^p < \infty ,
\]
(see (\ref{1.7}) in Theorem \ref {th4}).

\item
From the relation
$\|\cdot\|_{n,p} \leq \|\cdot\|_{m,p}$, $n\le m$,
it follows trivially
that
$\|\cdot\|_{-1-\eta,p}\le \|\cdot\|_p  $.
Thus, assumption 4 (b)  implies
\[
E\int_0^T ds \|f(s,\cdot,0)\|^p_{-1-\eta,p} < \infty ,
\]
(see again (\ref{1.7}) in Theorem \ref {th4}).

\item The above remarks show that the assumptions  of Theorem \ref{th4} and of
Theorem \ref{thm2.2} are fulfilled.
Hence, the existence of $u^W$ satisfying the weak formulation of
equation (\ref{1.4}) and $u^M$ satisfying the mild formulation is
assured.

Notice that
\[
\ch_p^{1-\eta} (T) \subset L_p(\Omega\times [0,T];L_p(\R^d)) .
\]
\end{enumerate}
\end{remark}


{\bf Proof of theorem \ref{thm4.3}}  Equation (\ref{4.4}) can be now written
as \begin{align*}
u^W(t,x)& = \int_{\R^d} dy \,u_0(y) G(t,x; 0,y)\\
& + \int_0^t ds\,\int_{\R^d} f\big(s,y,u^W(s,y\big) G(t,x; s,y))\\
& + \int_0^t W^kds\,\int_{\R^d} g^k\big(s,y,u^W(s,y\big) G(t,x; s,y)).
\end{align*}

We next prove that

\begin{equation}
\label{4.5}
E
\left(
\int_0^T dt \|u^W(t,\cdot)-u^M(t,\cdot)\|^p_p
\right)
= 0 .
\end{equation}

Indeed, from the equations satisfied by $u^w$ and $u^M$, respectively, it follows that
\[
E
\left(
\int_0^T dt \|u^W(t,\cdot)-u^M(t,\cdot)\|^p_p
\right)
\leq C(R_1 + R_2) ,
\]
with
\begin{align*}
R_1& =E \int_0^T dt\int_{\R^d} dx \Big |\int_0^t ds\int_{\R^d} dy
\big(
f(s,y,u^W(t,y))-f(s,y,u^M(t,y))
\big)\\
&\quad \times G(t,x;s,y)\Big|^p,
\end{align*}
\begin{align*}
R_2&= E \int_0^T dt \int_{\R^d} dx \Big| \int_0^t W^k(ds)
\int_{\R^d} dy \big( g^k(s,y,u^W(t,y))-g^k(s,y,u^M(t,y)) \big)\\
&\quad\times G(t,x;s,y) \Big|^p .
\end{align*}
Following the arguments of the proof of Theorem \ref{thm2.2} and
by virtue of the Lipschitz assumptions on $f$ and $h$, we obtain
\[
R_1 + R_2 \leq C\int_0^T dt\int_0^t ds E \left(\int_{\R^d} dx
|u^W(s,x) - u^M(s,x)|^p\right).
\]
We conclude by Gronwall's lemma applied to the function
\begin{equation*}
\Psi(T) = E \Big(\int_0^T dt \|u^W(t,\cdot) - u^M(t,\cdot)\|^p_p ,\big)
\end{equation*}
for $T\geq 0$.


The conclusion of Theorem \ref{thm4.3} can be strengthened
assuming, for instance,
instead of 4(b), that
$f(t,x,0)=h(t,x,0)=0$ and $u_0\in\cap_{p\ge 2}L_p(\R^d)$.
Indeed, in this case it has been proved
that a.s.
\[
u^W \in \cc^{\gamma_1,\gamma_2} ([0,T]\times \R^d),
\]
with $\gamma_1 < \frac{1-\eta}{2}$,
$\gamma_2 < 1 - \eta$, and we can show that $\omega$-a.s.
$u^M$ owns the same property.
\medskip

{\bf Acknowledgements:} The first author wishes to thank the Institut 
de Matem\`atica, Universitat de Barcelona and the 
Centre de Recerca Matem\`atica
(Bellaterra, Spain) for their support and warm hospitality.
The second author expresses her thanks to the Univertit\`a  degli Studi di Padova for a kind invitation on July 2003.

\end{document}